\input amstex

\loadbold
\documentstyle{amsppt}

\magnification=1200
\input xy
\xyoption{all}

\topmatter  

\title Infinitesimal $K$-theory
\endtitle
\author by Guillermo Corti\~ nas\endauthor
\affil Departamento de Matem\'atica\\
Facultad de Cs. Exactas\\
    Universidad de La Plata\\
       Argentina\endaffil

\address Departamento de
Matem\' atica, Facultad de Ciencias Exactas, Calle 50 y 115
(1900) La Plata, Argentina.
\endaddress
\email willie\@mate.unlp.edu.ar\endemail
\leftheadtext{Guillermo Corti\~ nas}
\rightheadtext{Infinitesimal $K$-theory}

\endtopmatter  

\document
\redefine\lim{\operatornamewithlimits{lim}}
\define\normal{\vartriangleleft}
\define\coker{\operatorname{coker}}
\define\fib{\twoheadrightarrow}
\define\defor{\tilde\twoheadrightarrow}
\define\cof{\rightarrowtail}
\define\C{C}
\define\pointed{\operatorname{Sets}_*}
\define\pairs{\C\vartriangleright\C}
\define\s{S}
\define\PC{\operatorname{Pro-}C}
\define\PS{\operatorname{Pro-}\s}
\define\LF{LF}
\define\eLF{\operatorname\LF}
\define\erF{\underline\eLF}
\define\totF{\underline\erF}
\define\LX{L\bold{X}}
\define\eLX{\operatorname\LX}
\define\erX{\underline\eLX}
\define\totX{\underline\erX}
\define\tot{\underset{=}\to{\operatorname{L}}}

\define\holi{\operatornamewithlimits{holim}}
\define\hoco{\operatornamewithlimits{hocolim}}
\define\ssets{\operatornamewithlimits{SSets}}
\define\ho{\operatorname{Ho}}
\define\jose{\operatorname{HoSSets}}
\define\hofiber{\operatorname{hofiber}}
\define\rat{\Bbb Q}
\define\loc{\C [{De\negthinspace f}^{-1}]}

\subhead 0. Introduction\endsubhead
\bigskip
In this paper we study the fiber $F$ of the rational Jones-Goodwillie 
character 
$$
F:=\hofiber(ch:K^\rat(A)@>>>HN^\rat(A))
$$
going from $K$-theory to negative cyclic homology of associative rings.
We describe this fiber $F$ in terms of sheaf cohomology. We prove
that, for $n\ge 1$, there is an isomorphism (Th. 6.2):
$$
\pi_n(F)\cong H^{-n}_{inf}(A,K^\rat)\tag{0}
$$
between the homotopy of the fiber and the hypercohomology groups
of $K^\rat$ on a non-commutative version of Grothendieck's infinitesimal
site ([Dix]). We fall short of expressing $F$ as the $K$-theory of 
a category. However we construct a natural map (7.4-7.6):
$$
K_*^\rat (Free^{inf})@>>>H^*_{inf}(A, K^\rat(Free))\cong\pi_*(F)\tag{1}
$$
between the $K$-theory of locally free modules on the non-commutative
infinitesimal site and hypercohomology $K$-theory of free modules.
These $K$-theory
results are obtained as a very particular case of a general
construction. The input of this construction consists of a category
$C$ equipped with a suitable notion of infinitesimal deformation,
(e.g. $C$=rings, deformation=surjection with nilpotent kernel) and
a functor $X:C@>>>$Spaces. The output is a natural map:
$$
H_{inf}(-,X)@>>>X
$$
where the source is generalized sheaf cohomology and maps deformations
into weak equivalences. Delooping the fiber of this map we obtain
a character $c^\tau:X@>>>\tau X$ which induces an isomorphism
of relative groups:
$$
c^\tau:\pi_*X(f)\cong\pi_*\tau X(f)\tag{2}
$$
for any given deformation $f$. We show that, for $C=$ associative rings
and $X=K^\rat$, the map $c^{\tau}_{n}:K^\rat_n(A)@>>>\tau_nK(A)$ is
the Jones-Goodwillie character ($n\ge 1$), and \thetag{2} is 
Goodwillie's isomorphism [G]. For 
general $C$ and
$KM=$group completion of a permutative category $M$, we obtain a map
$KM^{inf}@>>>H_{inf}(-,KM)$ of which \thetag{1} is a particular case. The 
hypercohomology construction is also interesting
for other functors. For instance if $C$ is the category of associative
algebras over a field of characteristic zero, then for the Cuntz-Quillen
de Rham supercomplex:
$$
H_{inf}(A,X_{dR})\cong HP_*(A)\tag{3}
$$
For $C=$commutative algebras and $X=$the commutative de Rham
complex $\Omega$:
$$
H^*_{inf}(A,\Omega)=H_{inf}(SpecA,O)\tag{4}
$$
The right hand side of \thetag{4} is the infinitesimal cohomology of
the structure sheaf in the sense of Grothendieck [Dix]. For some
functors $X$ such as the de Rham complexes \thetag{3} and \thetag{4},
 there is another interpretation for $H_{inf}(-,X)$. This is the derived
functor analogy of Cuntz-Quillen ([CQ2], [C1]). For associative
algebras over a field it consists of the following. Given an
algebra $A$ choose
a presentation $A=R/I$ as a quotient of a quasi-free algebra,
and put
$$
LX(A)=\holi_{\Bbb N^{op}}(n\mapsto X(R/I^n))\tag{5}
$$
We formalize and generalize this construction. First we show that a 
category with deformations $C$ admits
a structure akin to that of a closed model category which permits calculation
of the localization $C[{De\negthinspace f}^{-1}]$ as the homotopy
category of cofibrant objects. (For associative algebras over a field,
cofibrant=quasi-free; for commutative algebras, cofibrant=smooth). Then
we show that, under certain conditions, hypercohomology can
be computed as the left derived functor with respect to the
localization above (Prop. 5.5):
$$
H_{inf}(A,X)\cong LX\tag{6}
$$
For categories of rings the isomorphism above occurs whenever 
$X$ satisfies a power series version of the Poincar\'e lemma (Th. 3.5).
Such is notably the case of the de Rham complexes \thetag{3}, \thetag{4} in 
their respective categories of definition (Th. 4.0). We also show
that $X=BE^+\otimes\rat$, the rational plus construction of the elementary 
group satisfies \thetag{6} (the general linear group does not, as
there are problems with $K_1$).  Computationwise, this 
means that, for $n\ge 2$,  $H^{-n}_{inf}(A,K^\rat)$ can be computed as in 
\thetag{5} (Ths. 4.3, 6.2). 
\smallskip
This paper continues the line of research I started in [C1]. The main
results of [C1] reappear here either in a more general form (sections
2-4) or with shorter, simpler proofs (section 4).
\smallskip
The rest of this paper is organized as follows. In section 1
we develop the language of categories with deformations. In section 2
we construct the localization $C@>>>\loc$ for a category with
deformations. A general criterion for the existence of derived
functors with respect to this localization is established in section 
3. In section 4 we show that both non-commutative de Rham
cohomology and the rational $+$-construction of the elementary
group meet this criterion, and
compute their derived functors. The sheaf theoretic approach
is developed in section 5 where the character $c^\tau$ mentioned
above is  constructed.
The isomorphism \thetag{0} is proved in section 6. Also in this section
we conjecture that $H^n(A,K^\rat)=0$ for positive $n$, and
show that this conjecture is related to finding a non commutative
analogue of Grothendieck's isomorphism \thetag{4}. The map \thetag{1} is 
constructed in section 7; a more concrete interpretation of $Free^{inf}$
as a category of integrable connections with as maps the gauge
transformations is discussed.
\bigskip
\subhead{1. Categories with Deformations}\endsubhead
\bigskip
\definition{Shrinking functors in categories of interest 1.0}
By a {\it category of interest} we mean a pointed category $C$
which is closed under finite limits and finite colimits. In particular
kernels and cokernels exist in $C$. We further
require that a weak version of the first Noether isomorphism holds in $C$. We shall
presently make this condition precise. First 
recall that a map in a category of interest is called {\it normal} if 
it is a kernel of its cokernel, and is a quotient map (or is conormal)
if it is a cokernel of its kernel. We require that if $i:I\hookrightarrow A$
and $j:J\hookrightarrow A$ are normal and $i$ factors as $jk$ then the
following sequence is exact:
$$
0@>>>J/I@>>> A/I @>>> A/J@>>>0\tag{NIT}
$$
i.e. the first nonzero map is a kernel of the second, which is a cokernel
of the first.
A category of interest has products and co-products. We use Cuntz' notations 
for co-products; we write $QA=A*A$ for 
the coproduct of $A$ with itself, $\mu:QA @>>> A$ for the folding map, 
$qA:=\ker\mu$ and $\partial_i$ for the natural inclusions ($i=0,1$). Recall 
a subobject $N\subset A$ of an object $A\in C$ is
an equivalence class of monomorphisms. The fact that a subobject is normal
is indicated by $N\vartriangleleft A$.
Consider the category
$C\vartriangleright C$ with as objects the pairs $(A,N)$ consisting
of an object $A\in\C$ and a normal subobject $N\vartriangleleft A$, and
as maps $(A,N)@>>>(A',N')$ the maps $A@>>>A'$ in $\C$ sending $N$ into
$N'$. Note $\pairs$ is equipped with a natural projection $\pi:\pairs@>>>C$.
A {\it shrinking} functor in a category of interest is a functor
$\pairs@>>>\pairs$ preserving $\pi$, $(A,N)\mapsto (A,s(A,N))$, $f\mapsto f$, 
which satisfies the following:
\smallskip
\item{$s1$} $s(A,N)\subset N$
\smallskip
\item{$s2$} $s(\frac{A}{s(A,N)},\frac{N}{s(A,N)})=0$
\medskip
\noindent
Notice that the inclusion in $s1$ is necessarily normal, because $s(A,N)\normal A$
is. We shall be especially concerned with maps $f:A@>>>B$ such that 
$s(A,\ker f)$=0. We remark that, in view of (NIT), condition s2
says that $A/s(A,N)@>>>A/N$ is such a map. We shall use the following
notation for the powers of the functor $s$; we shall write $s(A,N)^n$ to
mean $s(A,s(A,s(\dots,s(A,N)\dots))))$ ($n$ times).
\enddefinition
\bigskip
\remark{Remark 1.1} One can equivalently define shrinking functors
as functors on the category of normal monomorphisms. To see the
equivalence, proceed as follows. Start
with a shrinking functor $s$ in the sense of the definition above.
Given any concrete normal mono (or monic) $\alpha:N@>>>A$, choose a representative
$s(A,N)@>>>A$ of $s(A,\text{class of\ \ }\alpha)$; if $N=0$ is the fixed
$0$ object ($C$ is pointed), choose $s(A,0)=0$. The construction is functorial
because if $\alpha: N@>>>A$ and $\beta:M@>>>B$ are monos and $f:A@>>>B$ is
a map such that $f\alpha$ factors through $\beta$ then the factorization
is unique. Conversely, if we start with a shrinking functor defined on
concrete monos $\alpha$, then the class of $s(\alpha)$ depends only on the
class of $\alpha$, as is easily checked. Further these constructions are
inverse to each other. Precisely starting with an $s$ on $\pairs$, then
going to the associated functor on monos and coming back induces the identity,
and the reverse composition is naturally equivalent to the identity. Hence
we may --and do-- identify both notions. Notice also that any shrinking
functor defined for normal monos can be extended to a functor from the arrow
category to the category of normal monos, by 
$f:B@>>>A\mapsto s(A,\ker(\coker f))$.\endremark
\bigskip
\example{Main Examples 1.2} Note that $s(A,N)$ does not have to
depend on $N$; we may always set $s(,)=0$. More interestingly,
if $C$ is the category of (not necessarily unital) associative rings,
then the square two sided ideal $s(A,N):=N^2$ is a shrinking functor. A 
shrinking functor which depends on both variables is $s(A,N)=<[A,N]>$,
the ideal generated by all commutators. Another is $s(A,N)=<[A,N]>+N^2$.
Lie rings also admit at least two interesting shrinking functors; we may either
set $s(A,N)=[N,N]$ or $s(A,N)=[A,N]$. 
The same is true of groups;
simply substitute Lie brackets by commutators. The name of shrinking is meant 
to
convey the idea that as we iterate the functor, the result is smaller.
Of course this need not happen in particular cases, as
we may have $s(A,N)^n=N$ for all $n$, e.g. if $N=A$ is a unital
associative ring and $s(A,N)=N^2$. However the property that $s(A,N)^n=0$
for some $n$, i.e. the `shrinkability' of $N$ in $A$ is interesting, as
it gives, depending on the choice of $s(,)$, the notions of nilpotency,
solvability and filtered commutativity. All these examples are particular
cases of the following.
\endexample
\bigskip
\example{General Example 1.3} 
Let $C$ be a category of interest, and
let $s$ be a shrinking functor. Then $r(A,N)=(A/s(A,N),N/s(A,N))$
is a functor, and is equipped with a natural map $\epsilon:1@>>>r$. Note that
the
map $\epsilon(A,N)$ has as kernel the `diagonal' subobject $\Delta s(A,N)=(s(A,N),s(A,N))$.
The functor $r$ maps $\pairs$ into the full subcategory $E\subset\pairs$
of those pairs $(A,N)$ such that $s(A,N)=0$, and $\epsilon$ is an isomorphism
on objects of $E$. It is not hard to check that $\pairs$ is a category of
interest if $C$ is one. Note also that $E$ is closed under kernels and cokernels.
Conversely, let $C$ be a category of interest, and let $E\subset\pairs$ be a full subcategory,
closed under kernels and cokernels.
Suppose a functor $r:\pairs@>>>E$ is given, together with a natural map 
$\epsilon: 1@>>>r$. Suppose further that for all $(A,N)\in\pairs$ the map
 $\epsilon(A,N)$ is a quotient map, has a diagonal subobject as kernel
and that if $(A,N)$ happens to live in $E$, then $\epsilon(A,N)$ is an
isomorphism. Then it is not hard to check that 
$(A,N)\mapsto (A,\Delta^{-1}(\ker\epsilon(A,N))$ is a shrinking functor.
Note further that these constructions are mutually inverse. Thus a
shrinking functor is the same thing as the data $(C,E,r,\epsilon)$. For 
instance in the Lie ring examples above, the choice $s(A,N)=[N,N]$
comes from choosing as $E$ the category of pairs $(A,N)$ with 
$N$ abelian. The choice $s(A,N)=[A,N]$ comes from the subcategory
$E'\subset E$ where in addition we require the action to be trivial. The
group and associative ring examples are similarly obtained.
\endexample
\bigskip
\example{Pro-Example 1.4} Let $\C$ be a category; we write $\PC$ for
the category of countably indexed pro-objects. Thus --by [CQ3]-- every object
of $\PC$ is isomorphic to an inverse system of the form
$\{ A_n:n\in\Bbb N\}$. We regard $\C$ as the full subcategory of constant
pro-objects in $\PC$. By [AM], $\PC$ is a category of interest if $\C$ is.
If $\C$ is of interest and is equipped with a shrinking functor $s(A,N)$,
then it is possible to extend $s(,)$ to all of $\PC$ in such a way that
if $\{ A_n:n\in\Bbb N\}$ is an inverse system with structure maps
$\sigma:A_*@>>> A_{*-1}$ and $N_*\vartriangleleft A_*$ is a collection
of normal subobjects such that $\sigma N_{n}\subset N_{n-1}$ then
$$
s(A,N)=\{s(A_n,N_n):n\in \Bbb N\}\tag7
$$
In fact since every representative of a normal subobject in $\PC$ is isomorphic 
--in the arrow category-- to an inverse system of normal sub-objects as
above, the formula \thetag{7} is almost the definition of a functor; 
there are however problems with the many choices involved. To get an 
unambigous definition proceed as follows. First extend the defintion of 
$s(A,N)$ from normal subobjects to arbitrary maps $f:B@>>>A$ as in 1.1.
The pro-extension of $(A,f)\mapsto s(A,f)$ is a well-defined functor, and
sends maps which are inverse systems of monomorphisms into inverse
systems of monomorphisms, whence it sends every  monomorphism
to a normal monomorphism, because every monomorphism is isomorphic
to one of such form in the arrow category. In particular we obtain a functor
mapping pairs 
$(A,N)$ of a pro-object $A$ and a sub-pro-object $N$ to a sub-pro-object 
$s(A,N)\subset N$; one checks that this functor satisfies 
s1 and s2 above. We write 
$s(A,N)^\infty:=\operatornamewithlimits{lim}_n s(A,N)^n$ for the inverse
limit; this, as well as the limit of any inverse system indexed by the natural
numbers exists in $\PC$ ([AM]). For example if $A$ and $N$ are as in (7) 
above, we have:
$$
s(A,N)^\infty:=\{ s(A_n,N_n)^n: n\in\Bbb N\}\tag8
$$
Note that $s(,)^\infty$ is idempotent as a functor 
$\PC\vartriangleright\PC@>>>\PC\vartriangleright\PC$; one checks further
it is again a shrinking functor.
\endexample
\bigskip
\bigskip
\subhead{Homotopy 1.5}\endsubhead Let $\C$ be a category with a shrinking functor 
and let $A\in\C$. We say that two maps $f,g\in\C (A,B)$ are {\it congruent},
--and write $f\equiv g$-- if $f*g s(QA,qA)^n=0$ for some $n\ge 1$. Thus if 
$s(,)$ is idempotent, --e.g. as in (8)-- then $f\equiv g$ iff $f*g$ 
factors through 
$QA/s(QA,qA)$. In general, $f\equiv g$ iff $f*g$ factors through
$$
CylA:=QA/s(QA,qA)^\infty
$$ 
in $\PC$. One checks that $\equiv$ is a reflexive and symmetric relation, 
and that it is compatible with composition on both sides in the restricted sense
that $f_0\equiv f_1$ $\Rightarrow$ $gf_0\equiv gf_1$ and $f_0h\equiv f_1h$
(whenever composition makes sense). It follows that the equivalence relation $\sim$ 
generated by $\equiv$ is compatible with composition in the ample sense that
$f_0\sim f_1$ and $g_0\sim g_1$ imply $f_0g_0\sim f_1g_1$ . We say that $f$ and $g$ are {\it homotopic} if $f\sim g$. In the particular
case when $s$ is idempotent, $CylA$ is already an object of $C$. Such is
the case of $\PC$ with either of the shrinking functors (7) and (8).
Both shrinking functors induce the same cylinder and therefore the same
homotopy relation.
Note also that, as coproducts in $\PC$ are indexwise 
--i.e. $Q\{A_i:i\in I\}=\{QA_i:i\in I\}$--
the pro-extension of the functor $Cyl$, 
$Cyl(\{A_i:i\in I\}=\{QA_i/s(QA_i,qA_i)^n: (i,n)\in I\times N\}$
is precisely the natural cylinder associated to the shrinking functor of $\PC$.
\bigskip
\proclaim{Lemma 1.5.1} Let 
$$
\xymatrix{&B\ar [d]^{p}\\
          R\ar@<1ex>[ur]^{f_0}
          \ar [ur]_{f_1}\ar[r]_f &A}
$$
be a commutative diagram in $\PC$. Suppose $s(B,\ker p)^\infty=0$.
Then $f_0\equiv f_1$.
\endproclaim
\demo{Proof} Consider the sum map $h=f_0*f_1:QR@>>>B$. We have 
$ph\partial_i=f$, (i=0,1); hence $p h=f*f=f\mu$. 
Thus $h$ maps $qR=\ker\mu\subset\ker f\mu$ into $\ker p$, and 
$s(QR,qR)^\infty$ into $s(B,\ker p)^\infty=0$. Therefore $h$ induces a homotopy
$CylR@>>>B$.\qed
\enddemo
\bigskip
\proclaim{Corollary 1.5.2} Let $p$ be as in the lemma above. Assume
further that $p$ admits a right inverse; then $p$ is a homotopy equivalence.
Precisely if $ps=1_A$ then $sp\equiv 1_B$.
\endproclaim
\demo{Proof} Apply the lemma with $R=B$, $f_0=1_B$, $f_1=sp$ and $f=p$.\qed
\enddemo
\bigskip
\definition{Fibrations 1.6} A class of {\it fibrations} in a category of
interest is a class $Fib$ of maps which  
satisfies the following variation of the dual of Waldhausen's
axioms for a class of cofibrations [Wa p.320]:
\roster
\item{Fib1} $Fib$ contains all isomorphisms.
\smallskip
\item{Fib2} $Fib$ is closed under composition and base change by arbitrary
maps. 
\smallskip
\item{Fib3} The folding map $A^{*n}@>>>A$ defined as the identity in each 
summand is a fibration ($A\in C$, $n\in\Bbb N$).
\endroster
As opposed to Waldhausen's fibrations, ours do not necessarily include
the map $A@>>>0$. On the other hand (Fib3) is not required in [Wa].
Fibrations shall be denoted by a double headed arrow $\fib$. 
If $\C$ and $Fib$ are as above, then we say that an object $A\in\C$ is
{\it rel-projective} if it has the left lifting property of [Q1] with
respect to $Fib$. We say that $\C$ (or rather $(\C,Fib)$) has
{\it enough rel-projectives} if every object $A\in\C$ is
the target of a fibration $P_A\fib A$ with $P_A$ rel-projective.
\enddefinition
\bigskip
\example{Underlying example 1.7} In most of the examples considered in 
this paper,
fibrations are those maps which admit a right inverse in an 
{\it underlying category}. By the latter we understand a fixed category 
of interest $\s$ together with a faithful embedding $\C\subset\s$ which 
preserves inverse limits. Hence we can
equip $\PC$ with the class of fibrations consisting of those maps having
a right inverse in $\PS$. In all the examples 1.1 we may take 
$S=\operatorname{Sets}_*$, the category of pointed sets; in the
ring examples, we may also choose $S=$Abelian Groups. If instead of rings
we look at algebras over some ground ring $k$, $S=k-\operatorname{Mod}$
is another natural choice.
If, as is the case in these examples, the embedding has a left adjoint, 
then $\C$ has
sufficient relatively projectives. For if $\perp :\C @>>> \C$ is the
associated cotriple, then the co-unit map $\perp A\fib A$ is a fibration
with rel-projective source. Note that all this structure is preserved by 
the pro-category. Indeed the pro-extension of a faithful functor is 
faithful because both filtrant direct and inverse limits preserve 
injections, and if $L$ is left adjoint to the inclusion $\C\subset\s$,
then its pro-extension if left adjoint to $\PC\subset\PS$ (cf. [AM]). We 
remark that this class of fibrations has the following extra property:
$$
fg\in \text{Fib}\implies f\in\text{Fib}\tag{Fib4}
$$
\endexample
\bigskip
\subhead{Categories with deformations 1.8}\endsubhead
By a {\it category with {\rm (infinitesimal, iterative)} deformations 
({\rm or} thickenings)} we shall understand a category of 
interest $C$ together with a shrinking functor and a class of fibrations.
We require further that the following axiom be satisfied:
\medskip
\noindent\thetag{Def} If $B\fib A$ is a fibration with kernel $I$ then 
$B/s(B,I)\fib A$ as well as each of the
maps $B/s(B,I)^{n+1}\fib B/s(B,I)^n$ ($n\in \Bbb N$)
is a fibration.
\medskip
\noindent Note that if $s$ happens to be idempotent only the first condition 
is relevant. For example if fibrations are as in 1.7, then $\PC$ satifies
\thetag{Def} (by \thetag{Fib4}); if further every quotient map in $C$ is split in the underlying
category, then $C$ satisfies \thetag{Def} also.
If $C$ is any category with deformations, and $A,B\in C$ are
objects, then by a {\it deformation} of $A$ by $B$ we shall mean a fibration
$f:B\fib A$ such that if $N=\ker f$ then $s(A,N)^n=0$
for some $n\ge 1$. Thus if $s(,)$ happens to be idempotent, the latter
condition simply means that $s(A,N)=0$. In general we have $s(A,N)^\infty=0$.
For example, it follows from NIT, s2 and Def that the map $B/s(B,\ker f)^n\defor A$ induced
by a fibration $f:B\fib A$ is always a deformation. In particular
$$
QA/s(QA,qA)^n\defor A\tag9
$$
is a deformation by Fib3.
\bigskip
\subhead{Cofibrancy 1.8.1}\endsubhead 
A map $A@>>>B\in C$ is called a {\it cofibration} if it has the left lifting 
property (in the sense of [Q1]) with respect to deformations; an object $A$
is called {\it cofibrant}
if $0\cof A$ is a cofibration. For example relatively projective
objects are always cofibrant; if $s=0$ they are all the cofibrants.
Here are some concrete examples. We fix Fib=all
surjections in all cases, except as noted.
If $C$ is the category of commutative algebras
over a field, and $s(A,I)=I^2$,
then cofibrancy is the same as smoothness;
for associative algebras, cofibrant=quasi-free in the sense
of Cuntz-Quillen. For commutative
algebras over a ring which is not a field, smooth=cofibrant+flat. If
$C$ is the category of groups, and $s(G,N)=[G,N]$ is the relative commutator
subgroup, 
then $G$ is cofibrant iff $H^2(G,M)=0$ for every trivial
$G$-module $M$. Equivalently, by the universal coefficient theorem,
$G$ is cofibrant iff $G_{ab}$ is free abelian and $H_2(G,\Bbb Z)=0$. In
particular superperfect groups (i.e. universal central extensions of perfect
groups) are cofibrant in this setting. Similar considerations
can be made regarding the category of Lie algebras over a field with
$s(L,N)=[L,N]$. As noted above, in the case of groups we may also
take $s(G,N)=[N,N]$; here cofibrant groups are those for which
$H^2(G,M)=0$ for all $G$-modules $M$. Clearly free groups satisfy this;
the converse is a theorem of Stallings and Swan (cf. [Co]). Thus the
cofibrants are just the free groups in this case.
One can take any of these
examples and pass to the pro-category, extending $s$ as in (8). Then
one can either take as fibrations the maps which are split in Pro-Sets,
or all effective epimorphisms. The latter choice gives more deformations and
therefore less cofibrants (cf. [CQ3]). 
\smallskip
The following elementary
lemma shall be of use in what follows.
\bigskip
\proclaim{Lemma 1.9} Let $C$ be a category with deformations. Then:
\smallskip
\item{1)} The class of deformations contains all isomorphisms and is 
closed under composition and under base change by arbitrary maps.
\smallskip
\item{2)} The class
of cofibrations contains all isomorphisms and is closed under 
composition and cobase change by arbitrary maps. 
\smallskip
\item{3)} An object $A$ is cofibrant iff every deformation $\pi:B\defor A$
is split; i.e. there exists $s:A@>>>B$ such that $\pi s=1$.
\endproclaim
\demo{Proof} The first assertion of 1) follows from Fib1, the fact that
isomorphisms have trivial kernel and the fact that --by $s1$-- $s(A,0)=0$ 
($A\in C$).
Let $f:A\defor B$ and $g:B\defor C$ be deformations.
Write $I=\ker f$ and  $J=\ker g$, and choose $n$ and $m$ such that 
$s(A,I)^n=0$ and $s(B,J)^m=0$. Then $gf$ is a fibration and $K=\ker gf$ 
is the pullback of $J\subset B$ along $f$.
Hence $s(A,K)^m$ maps to zero in $B$ and therefore is contained in $I$.
Thus $s(A,K)^{m+n}=0$, and $gf$ is a deformation, proving the second 
assertion. Let $\pi:X\defor Y$ be a deformation and let $f:A@>>>Y$ be any
map. Write $P$ for the pullback and $\hat{\pi}$ and $\hat{f}$ for the
induced maps. Then $\hat{\pi}$ is a fibration by Fib2. Further,  $\hat{f}$ 
maps $K=\ker\hat{\pi}$ isomorphically onto
$I=\ker\pi$ and $s(P,K)^n$ monomorphically into $s(X,I)^n$ which is
zero for $n>>0$. Therefore $\hat{\pi}$ is a deformation. Thus 1) is proven.
The dual statements of 2) 
are immediate from the definition of cofibration. The third assertion
follows from the fact that deformations are closed under base change.\qed
\enddemo
\proclaim{Corollary 1.10} (Homotopy Extension) If $R$ is a cofibrant object 
in $C$ and $X\defor Y$ $\in C$ is
a deformation then the dotted arrow in the commutative diagram below
exists in $\PC$:
$$
\xymatrix{R\ar@{>->}[d]_{\partial_0}\ar[r] & X\ar@{>>}[d]\\
          CylR\ar[r]\ar@{-->}[ur]& Y}
$$
\endproclaim
\demo{Proof} From the lemma the map $\partial_0:R\cof QR$ is a cofibration.
Hence there exists a map $h:QR@>>>X$ making the obvious diagram commute. As
$X\defor Y$ is a deformation, the map $h$ must kill some power $s(QR,qR)^n$;
the induced map is a pro-map $CylR@>>>X$ with the desired property.\qed
\enddemo
\bigskip  
\subhead{Cofibrant models 1.11}\endsubhead
We say that $C$ has 
{\it enough cofibrant objects} if for each object $A\in C$ there exists a 
fibration $RA\fib A$ for
some cofibrant object $RA$. Such is the case for example if $\C$ has
enough rel-projectives.
By a {\it model} of an object $A\in C$ we understand a deformation 
$R\defor A$ where $R$ is cofibrant; we say that $C$ has 
{\it enough models} if each object has a model. Such is the case for example 
if $\C$ has sufficient cofibrant objects and $s(,)$
is idempotent. Indeed if $\pi^A:RA\fib A$ is a fibration with $RA$ cofibrant
and $I:=\ker \pi^A$, then the induced map ${\pi'}^A:R'A:=RA/s(RA,I)\defor A$, 
which is a deformation by s2 and Def, is a model. To see this, we must show that
the dotted arrow out of $R'A$ in the diagram below exists for every deformation 
$p:X\defor Y$:
$$
\xymatrix{RA\ar[d]\ar@{-->}[r] &X\ar@{>>}[d]\\
          R'A\ar@{-->}[ur]\ar[r] &Y} \tag10
$$
But the arrow out of $RA$ exists by cofibrancy of the latter, and it
maps $s(RA,I)$ into $\ker p$ and $s(RA,I)=s(RA,s(RA,I))$ into $s(X,\ker p)=0$.
Hence $RA@>>>X$ factors through $R'A$ as claimed. The same argument
shows that if $R$ is any cofibrant object and $I\vartriangleleft R$ is
any normal suboject, then $R/s(R,I)$ is cofibrant. In the
general case, i.e.
when the shrinking is not idempotent, one needs to pass to the pro-category
to find sufficient models. 

\smallskip
The idea of cofibrant models is analogous to that of projective
resolutions in abelian categories such as modules over a ring,
and to simplicial resolutions in non-abelian ones such as Rings.
As with resolutions, we can define derived functors from cofibrant
models. To formalize this we need to have a derived category; this
is constructed in the next section.
\bigskip
\subhead{2. The derived category}\endsubhead

Throughout this section we work under the following:
\smallskip
{\bf STANDING ASSUMPTION:} $C$ is a category with deformations
having 

\noindent sufficient cofibrant objects.
\smallskip
The main purpose of this section is to prove the theorem below.
\bigskip
\proclaim{Theorem 2.0} Let $\C$ be a category with deformations; assume
$\C$ has sufficient cofibrant objects. Then the localization 
$\C [{De\negthinspace f}^{-1}]$ of $\C$ at the class of deformations
exists, and is equivalent to the category whose objects are all fibrations
$R\fib A$ with $R$ cofibrant, ($A\in C$), and where a map 
$p_1:R\fib A @>>> p_2:S\fib B$ is the homotopy class 
of a map of pro-objects 
$R/s(R,\ker p_1)^\infty @>>> S/s(S,\ker p_2)^\infty$ in $\PC$. If furthermore, 
$s(,)$
is idempotent, then $\C [{De\negthinspace f}^{-1}]$ is also equivalent
to the homotopy category of cofibrant objects.
\endproclaim
\smallskip
The proof of the theorem requires two lemmas.
 
\proclaim{Lemma 2.1} Let $f:A@>>>B$ be a map in $C$ and let $\pi^A:RA\fib A$
and $RB\fib B$ be fibrations with $RA$ and $RB$ cofibrant. Write
$R'A:=RA/s(RA,\ker\pi^A)^\infty$, $R'B:=RB/s(RB,\ker\pi^B)^\infty$
and ${\pi'}^A:R'A\defor A$ and ${\pi'}^B$ for the induced maps. Then:
\roster
\item There exists a map $\hat{f}:R'A@>>>R'B$ lifting $f$; i.e., 
such that ${\pi'}^B\circ\hat{f}=f\circ {\pi'}^A $.
\smallskip
\item Any two liftings of $f$ as in {\rm (1)} are congruent.
\endroster
\endproclaim
\demo{Proof} By cofibrancy of $RA$ and Def, one can lift $f$ to
a map $f'$ from the constant pro-object $RA$ into $R'B$. This map
sends $IA:=\ker\pi^A$ to $\ker{\pi'}^B$ and $s(RA,IA)^\infty$ to
zero. Thus $f'$ factors through a map $\hat{f}$, proving (1). The second
statement is immediate from 1.5.1.\qed
\enddemo
\medskip
\remark{Remark 2.1.1} Note that we do not affirm nor negate that 
$R'A$ is cofibrant. Such an assertion does not make sense as we
have not equipped $\PC$ with any particular class of fibrations.
However the same argument proves that the statement obtained by replacing
 $R'A@>>>A$ and $R'B@>>>B$ by arbitrary deformations with cofibrant source
holds. The point is that, unless we are in the idempotent case (cf. 1.11),
the standing assumption does not imply the existence of sufficient models
in $C$.
\endremark
\bigskip
\proclaim{Lemma 2.2} Let $f:A\defor B$ be a deformation. Then any
lifting $\hat{f}$ as in {\rm 2.1-1)} above is a homotopy equivalence.
\endproclaim
\demo{Proof} We keep the notations of the proof of the previous lemma.
In view of 2.1-2) it suffices to show that there is
one lifting which is a homotopy equivalence. We shall construct such
a lifting as a composite $R'A@>>>R'_*P_*@>>>R'B$, and show that each of
the two components admits an inverse up to congruence. Write 
$R'_0A=R'_0P_0=A$,
$R'_0B=B$. Construct by induction a commutative diagram:
$$
\xymatrix{R'_{n+1}P_{n+1}\ar@{>>}[r] 
          & P_{n+1}\ar@{>>}[r]\ar@{>>}[d]& R'_{n+1}B\ar@{>>}[d]\\
          & R'_{n}P_{n}\ar@{>>}[r] & R'_nB}
$$
Here the inner diagram is a pullback, and  the two maps out of 
$P_{n+1}$  are
deformations by 1.9. Write $R'_*P_*$ for the inverse system
$n\mapsto R_nP_n$ just constructed.
The pro-map $R'_*P_*@>>>R'B$ is
the first component of the lifting as announced above. Next we construct
a congruence inverse for this map. By cofibrancy of $RB$ and induction, the 
fibration $RB\fib B$ lifts to a family of maps $RB@>>>R'_nP_n$ making the 
following diagram commute:
$$
\xymatrix{R'_{n+1}P_{n+1}\ar@{>>}[r]
          & P_{n+1}\ar@{>>}[r]\ar@{>>}[dd]& R'_{n+1}B\ar@{>>}[dd]\\
        &&&RB\ar[dl]
        \ar[ul] \ar@{.>}[dll]\ar@/^1pc/@{.>}[ull]\ar@/^3pc/@{.>}[ulll]\\
         & R'_nP_n\ar@{>>}[r] & R'_nB}
$$
The map $RB@>>>P_{n+1}$ exists by cartesianity and the map 
$RB@>>>R'_{n+1}P_{n+1}$
is a lifting of the latter along the deformation $R'_{n+1}P_{n+1}\defor P_{n+1}$.
Write $l_n:RB@>>>R'_nP_n$ for the map just constructed. Then $l_n$ maps
$s(RB,IB)^n=\ker(RB@>>>R'_nB)$ into $J_n:=\ker(R'_nP_n\defor R'_nB)$, whence
$l_n(s(RB,IB)^{nm})\subset s(R'_nP_n,J_n)^m$ for all $m$. Since
by construction $s(R'_nP_n,J_n)^m=0$ for $m>>0$, we have that 
$l_n(s(RB,IB)^m)=0$ for $m\ge\alpha_n$ where $\alpha_n$ is some integer
sufficiently large which depends on $n$. Choose these numbers in such
a way that $\alpha_n<\alpha_{n+1}$. Since by construction 
$l=\{l_n:RB@>>>R'_nP_n\}$ is a map of inverse systems $RB@>>>R'_*P_*$, 
the induced map $l'=\{R'_{\alpha_n}B@>>>R'_nP_n\}$ is a map of inverse
systems also, going from the inverse system 
$R'_\alpha B:n\mapsto  R'_{\alpha_n}B$ to the inverse
system $R'_*P_*$. 
 Because by construction the sequence $\{\alpha_n\}$ is strictly
increasing, it is cofinal, whence $l'$ is a map of pro-objects 
$R'B@>>>R'_*P_*$. Furthermore the composite 
$R'_\alpha B@>>>R'_*P_*@>>>R'B$ is the natural
projection, which represents the identity of the pro-object $R'B$.
Hence the map $l':R'B@>>>R'_*P_*$ is a right inverse of the map 
$R'_*P_*@>>>R'B$ in $\PC$. Thus by 1.5.2, the maps
$l$ and $R'_*P_*@>>>R'B$ are homotopy inverse. We have thus constructed
one of two announced maps and proved it is a homotopy equivalence.
Next we construct the second map. Use the
cofibrancy of $RA$ to lift the fibration $\pi^A$ along the deformation
$R_1P_1\defor A$ first to a map $RA@>>>R_1P_1$ and then by induction
to a pro-map $t:RA@>>>R'_*P_*$ covering the identity of $A$. By a similar
argument as above, it is not hard to see that $t$ induces a map $t':R'A@>>>R'_*P_*$
still covering the identity of $A$. It remains to show that $t'$ has a homotopy
inverse. By 1.5.1 it suffices to show that there is also a map in the
opposite direction $u:R'_*P_*@>>>R'A$ covering the identity of $A$. 
The latter map shall be constructed as a composite 
$u:R'_*P_*@>\tau>>R'_\beta P@>v>>R'A$. We define the map $v$ as follows.
First use cofibrancy of $RP_1$ and the fact that $R'_1A\defor A$ is
a deformation to lift the composite $RP_1@>>>R'_1P_1\defor A$ to a
map $v'_1:RP_1@>>>R'_1A$. Then because $v'_1$ covers the identity of $A$, 
and because $R'_1A\defor A$ is a deformation, the map $v'_1$ factors through
a map $v_1:R'_{\beta_1}P_1@>>>R'_1A$ for some $\beta_1>1$. To construct
$v_2$, first lift $RP_2@>>>R'_2P_2\defor R'_1P_1$ along the deformation
$R'_{\beta_1}P_1\defor R'_1P_1$ to a map $RP_2@>>>R'_{\beta_1}P_1$.
Next lift the composite of the latter map followed by $v_1$ to a map
$v'_2:RP_2@>>>R'_2A$, covering the identity of $R'_1A$. By the same
argument as above, there is an integer $\beta_2>\beta_1$ such
that the map $v'_2$ factors through a map $v_2:R'_{\beta_2}P_2@>>>R'_2A$.
Proceeding inductively, we get an inverse system $R_{\beta_*}P_*$
and a map of inverse systems
$v:R'_{\beta_*}P_*@>>>R'A$ covering the identity of $A$, as claimed.
Next we need to find a map $\tau:R'_*P_*@>>>R'_{\beta_*}P_*$ covering
the identity of $A$. As a preliminary step, note that we already have a 
map in the opposite direction. Indeed by the construction of $R'_{\beta_*}P_*$,
the projection maps $\theta_n:R'_{\beta_n}\defor R'_nP_n$ commute with
structure maps, and therefore assemble into a map of inverse systems;
furthermore each $\theta_n$ is a deformation, by (Def). 
Next we shall construct the map $\tau$ and in so doing we shall show it
is an isomorphism inverse to the map $\theta$ just considered. Because
$\theta$ covers the identity of $A$, it will follow that the same
is true of $\tau$. Now to the construction of $\tau$.
Consider $RP_{\beta_{1}}$; by cofibrancy, there is a map from the 
latter to $R'_{\beta_1}P_1$ lifting the projection 
$RP_{\beta_{1}}@>>>R'_{\beta_{1}}P_{\beta_{1}}\defor R_1P_1$ along $\theta_1$.
By $NIT$ and $s2$, the kernel $K_1$ of $\theta_1$ satisfies 
$s(R'_{\beta_1}P_1,K_1)^{\beta_1}=0$; hence the indicated lifting maps
$s(RP_{\beta_1},IP_{\beta_1})^{\beta_1}$ to zero and hence factors 
through $R'_{\beta_1}P_{\beta_1}$. Consider the map 
$\beta:\Bbb N@>>>\Bbb N$, $\beta(n)=\beta_n$; proceeding inductively,
 we get a map of inverse systems 
$\tau_n:R'_{{\beta}^{n+1}(1)}P_{{\beta}^{n+1}(1)}@>>>R'_{{\beta}^{n+1}(1)}P_{{\beta}^n(1)}$,
which represents a pro-map $\tau:R'_*P_*@>>>R_{\beta_*}P_*$.
One checks that both $\tau\theta$ and $\theta\tau$ are restriction maps,
and therefore represent identity maps.\qed
\enddemo

\medskip
\remark{Remark 2.2.1} Note that, in the proof above, 
as $A\mapsto RA$ is not a functor--nor is $RA\fib A$ a deformation--
we do not really have a pro-object $RP$. This difficulty, as well as that
remarked in 2.1.1, can be
avoided if one assumes that fibrations in $C$ are as in the underlying
example 1.7, and that the forgetful functor has an adjoint. For then
the cotriple $RA=\perp A$ is functorial and $R'A$ is cofibrant by \thetag{10}.
\endremark
\bigskip
\medskip
\demo{\bf Proof of Theorem 2.0} For each $A\in\C$, choose a fibration 
$\pi^A:RA\fib A$ with
$RA$ cofibrant; if $A$ is cofibrant already, choose $\pi^A=1$. Write 
$IA:=\ker\pi^A$. Let $\C'$ be the category whose objects are those
of $\C$ and where a map $f:A @>>> B$ is the homotopy class of a map of
pro-objects $R'A:=RA/s(RA,IA)^\infty @>>> R'B$. Define a functor 
$\gamma:\C @>>>\C'$ as follows. Set $\gamma A=A$ on objects, and for each
map $f:A @>>> B$ choose a lifting $\hat f:R'A @>>> R'B$, and define
$\gamma f:=[\hat f]$ as the homotopy class of $\hat f$. Then $\gamma$ is
a functor by 2.1, and maps deformations into isomorphisms by 2.2.
Next, we have to prove that any functor
$F:\C @>>> D$ which inverts deformations factors uniquely as 
$F=\hat F\gamma$.
Define $\hat F:\C' @>>>D$ by $\hat FA=FA$ on objects; to define $\hat F$ 
on arrows, proceed as follows. Given a homotopy class $\alpha\in\C'(A,B)$ 
choose a representative $f\in\PC(R'A,R'B)$ of $\alpha$, then choose a 
representative $f_{nm}:R'A_m:=RA/s(RA,IA)^m @>>> R'B_n$ of $f$ and
finally set $\hat F\alpha=F(R'B_n\defor B)Ff_{nm}F(R'A_m\defor A)^{-1}$.
One checks immediately --using \thetag{9}-- that $\hat F\alpha$ is 
independent of 
the choices made in its definition and that it is functorial.
Next we have to prove that $\hat F$ is unique. Suppose $G:\C'@>>> D$ is
another functor such that $G\gamma=F$. Then $GA=FA=\hat FA$, hence
$\hat F$ and $G$ agree on objects; similarly, if $\alpha\in\C'(A,B)$
is the class of a map $g:R'A @>>> R'B$ which admits a representative
$f:A @>>> B$, then $G\alpha=G\gamma f=\hat F\alpha$. In general,
if $f:R'A_m @>>>R'B_n$ represents $g$, then the 
following diagram commutes in $\PC$:
$$
\CD
R'A @<\hat{\pi}^A_m<< R'(R'A_m) @>\hat f>> R'(R'B_n)@>\hat{\pi}_n^B>>R'B\\
@VVV                     @VVV            @VVV                      @VVV\\
A @<<\pi_m^A<          R'A_m  @>>f>        R'B_n @>>\pi_n^B>    B\\
\endCD
$$
Hence the following diagram commutes up to homotopy:
$$
\CD
R'(RA_m) @>\hat f>> R'(R'B_n)\\
@V\hat{\pi}_m^AVV   @V\hat{\pi}_n^BVV\\
R'A @>>g> R'B\\
\endCD
$$
\bigskip
Thus for $\alpha=[g]$ we have $G\alpha=G(\gamma({\pi_n}^B)\gamma
f{\gamma({\pi_m}^A)}^{-1})=\hat F\alpha$, and $G=\hat F$, as claimed.
This finishes the proof of the first assertion of the theorem. Next,
let $\C"$ be the category of fibrations $R\fib A$ as in the theorem. One
checks that $A\mapsto\pi^A$, $[f]\mapsto [f]$ is a fully faithful functor
$F:C'@>>>C"$. I claim further that every object in $C"$ is isomorphic to
one in the image of $F$. Note the claim implies that $F$ is a category
equivalence as we have to prove. To prove the claim, let 
$p: R\fib A\in C"$ be an object. Then by cofibrancy of $RA$ the map 
$\pi^A$ can be lifted
to a map $\hat p:RA @>>> R'=R/s(R,\ker p)^\infty\in PC$, which passes
to the quotient to give a map $\hat{p'}:R'A@>>>R'\in PC$, whose homotopy
class is a map $\pi^A@>>>p\in C"$. Next observe that by 2.1-2), the map
$[\hat{p'}]$ is an isomorphism, and the claim is proved. To finish
proof of the theorem it only remains to prove the last assertion; the
proof is similar to that just given for the next to last statement
of the theorem, so we shall just sketch it. Assume
that $s$ is idempotent; by its very definition, $\gamma$ induces
a fully faithful embedding $\gamma'$ of the homotopy category of cofibrant 
objects in the localization of $C$. But by 1.11 and (Def) every object in
the localization is isomorphic to one in the image of $\gamma'$. This
completes the proof.\qed
\enddemo
\smallskip
\remark{Remark 2.3} The proof of 2.0 does not require $C$ to be
of interest. For example the same proof as above applies for unital 
rings.
\endremark
\bigskip
\subhead{3. Derived Functors and Poincar\'e Lemma}\endsubhead
\smallskip
\noindent Throughout this section, $C$ is a fixed category with deformations
satisfying the Standing Assumption of \S3. We recall the following
definition from [Q1, Ch.I, \S 4.1].
\bigskip
\definition{Definition 3.0} Let $\C$ be a 
category with deformations,
let $\gamma:\C @>>> \C [{De\negthinspace f}^{-1}]$ be the localization
functor of 2.0 above, and let $F:\C @>>> D$ be a functor. By the
{\it total left derived functor} of $F$ we mean 
a functor $\totF:\loc @>>> D$ together
with a natural transformation $\alpha:\totF\gamma @>>>F$ having the 
following universal property. Given any $G:\loc@>>> D$ and natural
transformation $\beta:G\gamma  @>>>F$ there is a unique natural 
transformation 
$\theta:G@>>>\totF$ such that the following diagram commutes:
$$
\xymatrix{G\gamma\ar [r]^\beta\ar@{.>} [d]^{\theta\gamma} &F\\
          \totF\gamma\ar[ur]_\alpha&\\}
\tag{11}
$$
\enddefinition
\bigskip
The following lemma establishes a criterion for the recognition
of derived functors.
\bigskip
\proclaim{Lemma 3.1} With the notations of the definition above,
Let $C$ be a category with deformations having sufficient cofibrant objects.
Assume the shrinking functor
is idempotent. Let $F:C@>>>D$ be a functor. Then the following conditions
are equivalent:
\item{i)} $F$ carries deformations between cofibrant objects into isomorphisms.
\smallskip
\item{ii)} The map $F(\mu:CylR\defor R)$ is an isomorphism for every cofibrant
object $R\in C$.
\smallskip
\item {iii)} $F$ carries homotopic maps between cofibrant objects into equal maps.
\smallskip
\item{iv)} The following construction is independent --up to isomorphism--of the choices made in
its definition. Given an object $A\in C$ choose
a cofibrant model $RA\defor A$, and set $LFA=FR$. Given a map $f:A@>>>B$
choose a lifting $\hat{f}:RA@>>>RB$ and set $LF(f)=F(\hat{f})$.
\smallskip
\item {v)} The derived functor exists and $\totF\gamma(A)=F(A)$ for
all cofibrant objects $A$.
\smallskip
Under the equivalent conditions above, 
the construction $LF$ of {\rm (iv)} is functorial and $LF=\totF\gamma$.
\endproclaim
\demo{Proof} That i)$\implies$ii)$\iff$iii) is clear from the definition
of homotopy; ii)$\implies$i)
follows from 1.9-3) and 1.5.1.; iii)$\implies$iv) is immediate from 2.1.1.
If i) does not hold, then there is a deformation $f:R\defor S$ with
both $R$ and $S$ cofibrant such that $F(f)$ is not an isomorphism.
Hence the choices $RS=R$ and $RS=S$ lead to two distinct values of
$LF(S)$. Thus iv)$\implies$i). The last assertion of the lemma
implies that iv)$\iff$v). 
Note that, once $LF$ is assumed to be well-defined, it is automatically 
functorial,
and --by 2.2-- maps deformations into isomorphisms. Hence --by 2.0-- it
factors uniquely as $LF=G\gamma$. The proof that $G=\totF$ is essentially
the same as the proof of [Q1, I.4, Prop. 1]; details 
are left to the reader. \qed
\enddemo
\bigskip
\remark{Functors to Simplicial Sets 3.2} 
Suppose a functor $X:C@>>>\operatorname{SSets}$ is given. Let $A\in C$ and
let $RA$ and $IA$ be as in the previous section. Consider
the homotopy limit:
$$
LX(A)=\holi_{\Bbb N^{op}}(n\mapsto X(RA/s(RA,IA)^n)\tag12
$$
Note that this construction is analogous to the one used to define the
hyper(co) homology of a functor from some abelian category to its chain complexes.
Simply think of $R'A@>>>A$ as a resolution and of $\holi$ as a total complex.
More formally the homotopy type of $LX$ is analogous to the composite of the derived functor
in the sense of the definition above and the localization functor.
Note that, for \thetag{12} to make sense, we need to have a true functor to 
SSets,
as opposed to a homotopy type. However we at least want that if $X@>\sim>>Y$
is a natural weak equivalence, then $LX@>\sim>>LY$. For this we need $X$ and
$Y$ to be fibrant (i.e. Kan) simplicial sets (cf.[BK, XI 5.6]). Now if 
$X$ is fibrant
then by the proof of the cofinality theorem of [BK,XI 9.2], the following 
assignment is
a functor $\PC@>>>\operatorname{HoSSets}$ extending \thetag{12}: 
$$
\bold{X}\{ A_i:i\in I\}=\text{homotopy type of\ \ }\holi_IXA_i\tag{13}
$$
Suppose that a deformation category structure is given in $\PC$, such
that $R'A$ is cofibrant and $R'A\defor A$ is a deformation; such is
the case e.g. if fibrations are as in example 1.7. Then if the
equivalent conditions of 3.1 hold for \thetag{13}, 
the homotopy type of \thetag{12} is a functor, --by 3.1-iv)-- and is
precisely the restriction of the derived functor of \thetag{13} to the
subcategory
$C\subset\PC$. We remark that this construction, however useful, is not
really the
derived functor of the homotopy type of $X:C@>>>\operatorname{SSets}$,
which need not exist. This is because by defintion 3.0 the universal
property \thetag{11} must hold for functors $G$ with values in 
$\jose$, while the universal property of the $\holi$ in \thetag{12}
requires functors to SSets. On the plus side, if we happen to have
a functorial choice of $RA\fib A$ then \thetag{12} yields not just a
homotopy type but a true sset. Such is notably the case when fibrations
are as in the underlying example 1.7 and the forgetful functor has
a left adjoint. This example has the advantage that one does not need
to check the conditions of 3.1 for all the cofibrant objects of $\PC$,
as the next lemma shows.
\endremark 
\bigskip
\proclaim{Lemma 3.3} Let $C$ be a category of interest with a shrinking 
functor. Assume $C$ is equipped with a faithful functor $C@>>>S$ into
another category of interest $S$ which has a left adjoint.
Write $\perp$ for the associated cotriple, $JA:=\ker(\perp A@>>>A)$, and
$UA=\perp A/s(\perp A,JA)^\infty$. Consider
$C$ and $\PC$ as categories with deformations, with fibrations defined
as in {\rm 1.7} and shrinking functor as in \thetag{8}. Let $X:C@>>>$Fibrant 
SSets be a functor. Then the functor \thetag{13} satisfies the equivalent 
conditions of {\rm 3.1} iff the functor $X$ satisfies the following:
\smallskip
\item{vi)} For every object $A\in C$ the map $X(\mu(UA):CylUA@>>>UA)$
is a weak equivalence.
\smallskip
\endproclaim

\demo{Proof} Note that, by 1.5.2, $vi)$ is logically weaker than 3.1-iii). 
Hence
it suffices to show that $vi)$ implies at least one --and then all-- 
of the equivalent conditions
of 3.1.  Assume vi) holds; we shall prove that 3.1-i) holds also. 
It follows from vi) and the fact that $\holi$ preserves
weak equivalences of fibrant ssets ([BK, XI5.6]) that if $A$ is any
pro-object then the map $\bold{X}(\mu(UA))$ is an isomorphism in
$\jose$. Thus for $A,B\in\PC$, 
$\bold{X}:\PC(UB,A)@>>>\jose(\bold{X}(UB),\bold{X}(A))$ sends
homotopic maps to equal maps. Thus $\bold{X}$ sends the following maps
into isomorphisms:
\smallskip
\item{-} all homotopy equivalences $UB@>>>UA$,
\smallskip
\item{-} those homotopy equivalences $UB@>>>A$ which admit a strict --not 
just homotopy-- right inverse.
\smallskip
Now let $p:B\defor A\in\PC$ be a deformation of cofibrant objects;
we must prove $\bold{X}(p)$ is an isomorphism. Consider the following
diagram:
$$
\xymatrix{UB\ar [d]^{\pi^B}\ar [r]^{Up}&UA\ar[d]_{\pi^A}\\
          B\ar@/^1pc/ [u]\ar [r]_p &A\ar@/_1pc/[u]}
$$
Here the arrows going up are right inverses of those going down; they
exist by cofibrancy of both $A$ and $B$. By 1.5.2 and 2.2, each of
$\pi^B$, $\pi^A$ and $Up$ is a homotopy equivalence, whence   
$\bold{X}$ maps each of them --and therefore also the map $p$-- 
to an isomorphism.\qed
\enddemo
\bigskip
\bigskip
\example{Associative Rings 3.4} Let $C$ be the category of 
--not necessarily unital--
associative rings, equipped with the shrinking functor $s(A,I)=I^2$ and
with as fibrations the surjective maps. Then the hypothesis of the lemma
above are satisfied for $S=\pointed$, the category of
pointed sets, and the forgetful functor $C@>>>S$. We shall show below
(Theorem 3.5) that condition vi) of 3.3 is equivalent to a Poincar\'e lemma 
for power series. The key observation needed to prove theorem 3.5 is that,
as we shall see presently, the pro-ring $CylUA$ is a power series pro-ring. 
To make this assertion precise, we need some notation. Given a pointed set $S$
and a ring $B$, we write $B\{ S\}$ for the ring of polynomials in the non 
commutative variables $S-\{ *\}$ (we identify $*$ with $0$), and 
$<S>\subset B\{ S\}$ for the two-sided ideal generated by the variables. 
We think of the power series on $S$ as a pro-ring; for each (pro-) ring
$B$ and pointed set $S$, we put:
$$
B\{\{S\}\}=B\{S\}/<S>^\infty\tag{14}
$$
The rest of this subsection will be devoted to proving that there is
a natural isomorphism:
$$
CylUA\cong UA\{\{A\}\}\tag{15}
$$
Our proof uses Lemma 3.7, which is proved separately below.
We need some more notation. Write $V$ for the free abelian group
on $A-\{ 0\}$, $TV$ for the tensor algebra (over $\Bbb Z$), and let $I$
be the kernel of the adjunction map $TV@>>>A$. Thus in the notation
of Lemma 3.3, $TV=\perp A$, $I=JA$, and $UA=TV/I^\infty$. The left 
hand side of \thetag{15} is a quotient
of $QTV$; precisely, identifying $TV$ with its image through the
inclusion $\partial_1:TV@>>>QTV=T(V\oplus V)$, $\partial_1(v)=(v,0)$, 
$v\in V$, we have 
$Cyl(TV/I^\infty)=QTV/{\Cal F}^\infty$, where ${\Cal F}^\infty$ is 
the pro-ideal  
${\Cal F}^1=<I>+<q(I)>+(qTV)\supset\dots\supset
{\Cal F}^n=<I^n>+<q(I^n)>+(qTV)^n\supset\dots$. Similarly, 
we have $UA\{\{A\}\}=\{ (TV/I^n){A}/<A>^n:n\in\Bbb N\}=QTV/{\Cal G}^\infty$,
 where ${\Cal G}^n:=<I^n>+<\partial_0(V)>^n$, and $\partial_0(v)=(0,v)$. 
Consider the isomorphism $\alpha:QTV=T(V\oplus V)@>>>QTV$, 
$\alpha(v,w):=(v+w,-w)$, $(v,w)\in V\oplus V$. We shall show that 
$\alpha$ maps the filtration $\Cal F$ into a filtration equivalent
to $\Cal G$, and thus induces an isomorphism as in \thetag{15}. For this
purpose we consider four new filtrations by ideals of $QTV$. We put
${\Cal F'}^n=(<I>+qTV)^n$, ${\Cal F"}^n=<I^n>+(qTV)^n$, 
${\Cal G'}^n=(<I>+<\partial_0(V)>)^n$, ${\Cal G"}^n=<I^n>+<\partial_0(V)>^n$.
One checks that the isomorphism $\alpha$ maps $\Cal F'$ isomorphically onto $\Cal G'$.
Hence it suffices to show that $\Cal F$ is equivalent to $\Cal F'$ and
that $\Cal G$ is equivalent to $\Cal G'$. But since
${\Cal F'}^n\supset\Cal F^n\supset{\Cal F"}^n$ and
${\Cal G'}^n\supset\Cal G^n\supset{\Cal G"}^n$ we are reduced to showing
that, for $N$ sufficiently large, we have ${\Cal F"}^n\supset{\Cal F'}^N$
and ${\Cal G"}^n\supset{\Cal G'}^N$; both these inclusions follow
from Lemma 3.7 below.
\endexample
\bigskip
\proclaim{Theorem 3.5} (Poincar\'e Functors) Regard the categories of 
associative rings and pro-rings as deformation categories as in {\rm 3.4} above. Let 
$X:$Rings$@>>>$Fibrant Ssets be a functor. Extend $X$ to a functor 
$\bold{X}:$pro-Rings$@>>>$Ho-Ssets as in \thetag{13}.
We have:
\item{1)} If any of the following
holds, then $\totX$ exists and $\totX\gamma(A)=\bold{X}(A)$ for all quasi-free
pro-rings $A$:
\smallskip
\item{i)} (Poincar\'e Lemma for pro-power series) For every ring
$A$ and every set $S$, the map 
$\bold{X}(A\{\{S\}\}\defor A)$ associated to the pro-power series \thetag{14}
is an isomorphism.
\smallskip
\item{ii)} (Poincar\'e Lemma for polynomials) 
$\bold{X}(A[t]@>>>A)$ is an isomorphism for all $A\in$Rings.
\smallskip
\item{iii)} Either of the above holds for all pro-rings of the form
$R/I^\infty$ where $R$ is a quasi-free ring and $I\normal R$ is an ideal.
\bigskip
\item{2)}  If $\totX$ exists and $\totX\gamma(A)=\bold{X}(A)$ for all 
quasi-free
pro-rings $A$, then i) above holds for all quasi-free pro-rings. In such
case we say that $X$ is a {\rm Poincar\'e} functor.
\endproclaim
\demo{Proof} If 1-i) holds then $\totX$ exists and has the desired property
by 3.3 and 3.4. If 1-ii) holds and 
$A=A_0\oplus A_1\oplus A_2\dots$ is a graded ring, then $X$ maps the projection
$A@>>>A_0$  to a weak equivalence. As $A\{\{ S\}\}$
is pro-graded, it follows that 1-ii)$\implies$1-i). If (iii) holds then
so does the condition of 3.3, as $CylUA$ is of the indicated form, cf. 
\thetag{15}. This proves 1). If
$\totX$ exists and is as in 2), then by 3.1-i), $\bold{X}$ must map all 
deformations of 
quasi-free pro-rings into isomorphisms. One checks that $A\{\{S\}\}$
is quasi-free if $A$ is, whence 2) follows.\qed
\enddemo
\bigskip
\remark{Remark 3.6} By essentially the same arguments as above, it
is not hard to see that a functor $F:$pro-Rings$@>>>$Any Category 
satisfies 3.1 iff it satisfies 3.5-2). 
\endremark
\bigskip
\proclaim{Lemma 3.7} Let $A\subset B$ be rings and let $\epsilon:B\rightarrow A$
be a homomorphism such that $\epsilon a=a, (a\in A)$. Set $I=\ker\epsilon$,
and let $J\subset A$ be an ideal. Consider the following filtration in
$B$:
$$
B\supset {\Cal F}^n=<J^n>+I^n
$$
Then there is an isomorphism:
$$
B/{\Cal F}^\infty\cong B/(<J>+I)^\infty
$$
\endproclaim
\smallskip
\demo{Proof} Let ${\Cal G}^n=<J>^n+I^n$. It is straightforward to check
that $(<J>+I)^{2n}\subset{\Cal G}^n$, whence $B/(<J>+I)^\infty\cong B/{\Cal G}^\infty$.
Thus we must prove that $B/{\Cal G}^\infty\cong B/{\Cal F}^\infty$. It is clear
that ${\Cal G}^n\supset {\Cal F}^n$. I claim that for $N=n^2+n-1$, we also have
${\Cal G}^N\subset {\Cal F}^n$. To prove the claim --and the lemma-- it
suffices to show that $<J>^N\subset{\Cal F}^n$. Every element of $<J>^N$ is
a sum of products of the form:
$$
(j_1+i_1)\dots(j_N+i_N)\qquad\qquad (j_r\in J, i_r\in I)
$$
If we distribute all parenthesis we get a sum in which all  those
terms not in $I^n$ have at most n-1 factors in $I$ and at least $n^2$ 
factors in $J$. We must show that at least $n$ of the latter
appear side by side, forming a string. Assume the contrary holds. Then every
string of $j$'s must be broken off before the $n$-th step by an $i$. The
minimum number of $i$'s that are necessary for this to happen occurs when
each sequence of $j$'s is broken off exactly after the $n-1$-th $j$ in a 
row. Since there are $n^2=(n+1)(n-1)+1$ $j$'s, the minimum number must be
$n+1$, which is a contradiction.\qed\enddemo
\bigskip
\subhead{4. The derived functors of rational $K$-theory and Cyclic Homology}\endsubhead
\bigskip
In this section we give examples of derived functors in the sense
of the previous section. The first one is that which motivated
this paper.
\bigskip
\proclaim{Theorem 4.0} Consider the Cuntz-Quillen supercomplex 
{\rm [CQ2 (1)]}
as a functor $X:\rat -$Algebras$@>>>$Supercomplexes. Then $X$ is Poincar\'e,
and its derived functor is represented by the homotopy type of the periodic
cyclic complex.
\endproclaim
\demo{Proof} Immediate from [CQ2,(9)], [CQ2, 8.1] and 3.5-6 above.\qed
\enddemo
\bigskip
\remark{Remark 4.1} One can also derive the commutative de Rham complex
in the category of commutative rational algebras. The resulting derived
functor is the infinitesimal cohomology of Grothendieck [Dix], also 
called algebraic de Rham ([H]) and crystalline ([FT]).
\endremark
\bigskip
\subhead{4.2. The derived functor of rational $K$-theory}\endsubhead
Recall Goodwillie's isomorphism:
$$
K_*^\rat(A,I)\cong HN_*(A\otimes\rat,I\otimes\rat)\tag{16}
$$ 
between the relative rational $K$-group of
a nilpotent ideal and its analogue in negative cyclic homology [G]. 
Using this isomorphism, we show in Theorem 4.3 below that $K^\rat$ is (almost)
a Poincar\'e functor and that its derived functor is essentially the
fiber of the Jones-Goodwillie character:
$$
ch_*\otimes\rat :K_*^\rat(A)@>>>HN_*(A)\otimes\rat\tag{17}
$$
In order to apply the framework developed in the previous section,
we need some preliminaries. First of all we need a fibrant functorial
model for $K$-theory. One such model is 
$K=\Bbb Z_\infty NGl$, the
Bousfield-Kan completion of the nerve of the general linear group. As is
well-known, this
is just a functorial plus construction. To get
rational $K$-theory we complete again: $K^\rat=\rat_\infty K=\rat_\infty NGl$.
Actually $\rat_\infty X$ is fibrant for every sset $X$, so we can 
choose any other --not necessarily fibrant-- functorial model for integral 
$K$-theory.
Next note that the character \thetag{17} is induced by a natural map of
fibrant ssets $ch:K^\rat@>>>SN^\rat$. For example in the plus construction
approach, the simplicial map $ch$ is constructed as follows. Start off with
the Hurewicz map $NGl@>>>\Bbb Z NGl$. Next consider the Dold-Kan functor
$S:$Chain Complexes of abelian groups $@>>>$S.abelian groups, and follow
the
Hurewicz map with the result of applying $S$ to the 
chain map $\Bbb Z NGl@>>>CN_{\ge 1}$ of [G]. We thus get a map
$NGl@>>>SCN_{\ge 1}$. Finally $\rat$-complete on both sides to obtain
a map $ch:K^\rat@>>>\rat_\infty SCN_{\ge 1}=:SN^\rat$. This is essentially
construction of [We], modulo geometric realization. Note for example
that the map $SCN_{\ge 1}\otimes\rat@>\sim>>SN^\rat$ is an equivalence
--by virtue of [BK, V 3.3]-- whence $\pi_*SN^\rat\cong HN_*\otimes\rat$.
Now the map $ch:K^\rat@>>>SN^\rat$ is defined for
unital rings. We extend it to non-unital rings in the usual manner,
i.e. by considering the unital ring $\tilde{A}=A\oplus\Bbb Z$ and taking the
homotopy fiber of the simplicial map associated to the projection
$\tilde{A}@>>>\Bbb Z$. One checks that if $I\normal A$ is a nilpotent
ideal in a not necessarily unital ring $A$, then $K_n^\rat(A,I)=K_n^\rat(\tilde{A},I)$
and $HN_n^\rat(A,I)=HN_n^\rat(\tilde{A},I)$ (as observed in [Wo, proof of 
Th. 1.1]) whence \thetag{16}
holds for not necessarily unital rings. Next extend both $K^\rat$ and
$SN^\rat$ to pro-rings as in \thetag{13}. We observe that, as $\holi$
preserves homotopy fibration sequences of fibrant ssets, the isomorphism
\thetag{16} holds for arbitrary deformations of pro-rings. In particular
for
the pro-power series \thetag{14}
we have $K_n^\rat(A\{\{S\}\},<S>)\cong HN_n(A\otimes\rat\{\{S\}\},<S>\otimes\rat)$.
Now if $A$ is quasi-free as a ring then $A\otimes\rat$ and $A\otimes\rat\{\{S\}\}$
are quasi-free as $\rat$-algebras. Hence $K_n^\rat(A\{\{S\}\},<S>)=0$
for $n\ge 2$, as $HN_n$ of a quasi-free algebra is zero in degrees $\ge 2$.
However $K_1^\rat(A\{\{S\}\},<S>)=\ker(HH_1(A\{\{S\}\}@>>>HH_1(A))\ne 0$
in general, whence $K^\rat$ is not Poincar\'e. To get a Poincar\'e functor
out of $K^\rat$ we just need to eliminate its first homotopy group. 
We do this by substituting the elementary group for $Gl$; i.e.
we consider 
$$
KE^\rat(A):=\rat_\infty NE(A)
$$
\bigskip

\proclaim{Theorem 4.3}{\rm (The derived functor of $K$-theory)}
\smallskip
The functor $A\mapsto K^\rat(A)$ is not Poincar\'e.
However, the functor $A\mapsto KE^\rat (A)$ above is, and therefore
it has a left derived functor $LKE^\rat$. Set $LK^\rat_n(A):=\pi_nLKE^\rat$; 
then:
\item {i)} There is an exact sequence:
$$
\split
\dots HN_{n+1}A@>>>LK_n^\rat(A) @>>>K_n^\rat(A)@>>>HN_n(A)@>>>\\
\dots HN_3(A) @>>>LK_2^\rat(A)@>>>K_2^\rat(A)@>>>HN_2(A)\\
\endsplit
$$
\item{ii)} If $A=R/I$ is a presentation of $A$ where $R\otimes\rat$
has Hochschild dimension $\le n-1$, --i.e. $HH^{n}(R\otimes\rat,-)=0$--, 
then $LK^\rat_{n}(A)=\pi_{n}(\operatornamewithlimits{lim}_m K^\rat(R/I^m))$.
\endproclaim
\demo{Proof} The first two assertions follow from the discussion above.
To prove i), consider the exact sequence of $K$-groups
associated with the deformation $\pi^A:UA\defor A$. Then 
$LK_n^\rat(A)=K_n^\rat(UA)$ ($n\ge 2$) and $K_n(\pi^A)\cong HN_n(\pi^A)$
($n\ge 1$) . On the other hand $HN_n(UA)=0$ for $n\ge 2$, and therefore 
$HN_n(\pi^A)\cong HN_{n+1}(A)$, for $n\ge 2$. This proves that the sequence 
is exact at $LK_2^\rat(A)$ and to the left. By the same argument as above, 
the natural map $HN_2(A)\hookrightarrow HN_1(\pi^A)$ is injective, whence 
$K_2^\rat(A)@>>>K^\rat_1(\pi^A)$ factors through $ch_2$. It follows that the
sequence is exact also at $K_2^\rat(A)$, completing the proof of i). If $A=R/I$
is a presentation as in ii), then $HN_m(R/I^\infty)\otimes\rat=0$ for $m\ge n$.
Hence by i), $LK_{n}(A)=LK_{n}(R/I^\infty)=K_{n}(R/I^\infty)$.
\qed\enddemo
\bigskip
\remark{Remark 4.4} For commutative algebras over $\rat$ part ii) of the 
theorem applies to presentations $A=R/I$ with $R$ smooth of Krull dimension 
$<n$.
Indeed the latter coincides with Hochschild's for smooth algebras. In 
general, for any $R$ as in ii), we have an exact sequence ([BK]):
$$
0@>>>\operatornamewithlimits{lim^1}_mK^\rat_{n}(R/I^m)@>>>LK^\rat_{n}(A)\fib\operatornamewithlimits{lim}_mK^\rat_{n}(R/I^m)@>>>0
$$
As the map $LK^\rat_n(A)=K^\rat_n(R/I^\infty)@>>>K^\rat_n(A)$ is induced by 
the projection 

\noindent $R/I^\infty\defor A$,
it factors through the lim term above. Hence the $\operatorname{lim}^1$ term 
is an obstruction for the surjectivity of the character 
$ch:K^\rat_{n+1}(A)@>>>HN_{n+1}(A)$. \endremark
\bigskip
\bigskip 
\subhead{5. Sheaf Theoretic Approach}\endsubhead

In this section we go back to the general setting of sections 1-3; 
we work in a fixed category $C$ with deformations. In the previous
section we have seen the virtues and pitfalls of the derived
functor construction of section 3. The pitfall being mainly that
$\tot F$ may not exist for a given functor $F$. In this section we
produce another object; the infinitesimal hypercohomology of a functor
to simplicial sets. This construction is closely related to the derived
functor construction, and has the advantage of being defined without
any hypothesis on the functor in question.
\bigskip
\subhead {5.0 Infinitesimal cohomology}\endsubhead
Given $A\in C$ consider the category
\noindent $inf(C/A)$ with as objects the deformations $B\defor A$ and as maps
$(B_0\defor A)@>>>(B_1\defor A)$ the maps $B_0@>>>B_1\in C$ making the 
diagram commute. Assume $C$ --whence also $inf(C/A)$-- is small. 
We regard $inf(C/A)^{op}$ as a site with the indiscrete topology; i.e. as coverings we take the 
families $\{B\cong B'\}$ consisting of a single isomorphism. Thus 
a sheaf on $inf(C/A)^{op}$ is just any covariant functor on $inf(C/A)$. 
Recall that if $G$ is a sheaf of abelian groups, then its cohomology groups are defined
as the right derived functors of its global sections 
$\operatornamewithlimits{lim}_{inf(C/A)}G$. 
By analogy, if $X$ is a sheaf of simplicial sets, we define its 
hypercohomology
as the right derived functor of its $\holi$. Precisely, for each $A$
the category
$SSets^{inf(C/A)}$ is a closed model category ([BK, proof of XI 8.1]), so we take the 
total right derived functor:
$$
H_{inf}(A,X)=\underset{=}\to{\operatorname{R}}\holi_{inf(C/A)}X
$$
Although this definition makes sense in general, we shall apply it
for ssets with the property that $\pi_nX$ is an abelian group
for all $n$. 
If $X$ is fibrant (i.e. if $X(B)$ is fibrant
for all $B\defor A\in inf(C/A)$) this homotopy type is calculated by
$H_{inf}(A,X)\cong\holi_{inf(C/A)}X$. If $X$ is any sheaf,
then $H_{inf}(A,X)=\holi_{inf(C/A)}X'$ where $X'$ is any
fibrant sheaf with a cofibration and weak equivalence 
$X\tilde{\rightarrowtail}X'$. 
For functors $X:C@>>>SSets$, this construction has properties in common 
with the derived
functor of section 3. For example $H_{inf}(A,X)$ always maps
deformations into isomorphisms. This follows from the cofinality
theorem for $\holi$, once one observes that if $f:B\defor A$ is
a deformation then $f_*:inf(C/B)@>>>inf(C/A)$ is left cofinal in the sense
of [BK, XI 9.1]. In turn, the cofinality of $f_*$ is immediate from the fact 
that, for each $\pi:E\defor A\in inf(C/A)$, 
the pullback $P\defor A$ of $\pi$ along $f$ is a final object of 
the over category $f/\pi$.
Another common feature between $L$ and $H_{inf}$ is that, if 
$X$ is fibrant, then 
 $H_{inf}(A,X)=\holi_{inf(C/A)}X@>>>\holi_{inf(C/A)}X(A)@>>>X(A)$
 is a natural map, and is a weak equivalence iff $X$ maps deformations into
 weak equivalences. A difference between $H_{inf}(A,X)$ and $LX$
 is that the first one exists independently of any property of $X$, appart
 from the smallness of $C$. Another is simply that $H_{inf}(A,X)$ does 
 not have the universal property of $LX$. Namely if $Y@>>>X$ is a map
 of fibrant sheaves, we have a commutative diagram:
 $$
 \CD
 H_{inf}(A,Y)@>>> H_{inf}(A,X)\\
 @VVV                    @VVV  \\
 Y(A)@>>>X(A)\\
\endCD
 $$
If $Y$ maps deformations into weak equivalences then the first vertical
map is a weak equivalence. This implies
that it has a --perhaps not natural-- homotopy inverse, which we can
compose with the first horizontal map to get a map $YA@>>>H_{inf}(A,X)$
making the diagram commute up to homotopy. To have a natural homotopy 
inverse of the first vertical map we need that $Y$ be cofibrant
as an object of $\ssets^C$, which is a very rare property. In fact the
latter is a model category, so one can replace any object 
by one that is cofibrant in this sense. This means that the lifting
exists in the homotopy category of $\ssets^C$. But we do not want
an object of $\ho(\ssets^C)$, we want a true functor, i.e. an object
of $\ho\ssets^C$. 
So in general $H_{inf}$ and $L$ are distinct.
A more precise comparison between these is given in 5.5 below. We use
the following indexing:
$$
H^n_{inf}(A,X):=\pi_{-n}H_{inf}(A,X)
$$

\bigskip
\proclaim{Lemma 5.1}(\v Cech pro-covering) Let $C$ be a --not necessarily small-- category with
deformations, and let $f:R\fib A$ be a fibration. Consider the n-fold
sum map $f^{*n}:R^{*n}\fib A$; write $q_n=\ker f^{*n}$. If $R$ is cofibrant, 
then the
functor $\check{f}:\Delta\times\Bbb N^{op}@>>>inf(C/A)$, 
$(n,m)\mapsto R^{*n+1}/s(R^{*n+1},q_{n+1})^m$
is left cofinal.
\endproclaim
\demo{Proof} We have to show that, for every object $\pi:B\defor A\in inf(C/A)$
--$B$ for short-- the category $\check{f}/B$ is null homotopic. By definition
an object of $\check{f}/B$ is a pair $(n,m)\in\Delta\times\Bbb N^{op}$
together with a map $\alpha:R^{*n+1}/s(R^{*n+1},q_{n+1})^m@>>>B\in inf(C/A)$.
 Let $I=ker\pi$, and
let $r\ge 1$ such that $s(B,I)^r=0$. Let $g:\Delta\times\Delta^r@>>>inf(C/A)$
be the restriction of $\check{f}$. We have a functor $\theta:\check{f}/B@>>>g/B$
given by $(n,m,\alpha)\mapsto(n,\min(m,r),\hat{\alpha})$, where $\hat{\alpha}$
is the map induced by passage to the quotient by $s(R^{*n+1},q_{n+1})^{min(m,r)}$.
There is also a natural faithful inclusion $\iota:g/B\subset\check{f}/B$.
We have $\theta\iota=1$; quotient by $s(-,-)^{min(-,r)}$ gives a natural
map $1@>>>\iota\theta$. Hence $\hat{f}/B$ is homotopy equivalent to
$g/B$. Next consider $h:\Delta@>>>inf(C/A)$, $n\mapsto g(n,r)$. Then
$\theta':(n,m,\alpha)\mapsto (n,\frac{R^{*n+1}}{s(R^{*n+1},q_{n+1})^r}\defor \frac{R^{*n+1}}{s(R^{*n+1},q_{n+1})^m}@>\alpha>> B)$
is a left inverse of the natural inclusion $\iota':h/B\subset g/B$,
and is equipped with a natural map $\iota'\theta'@>>>1$. Hence $g/B$ is
weakly equivalent to $h/B$. Now the deformation 
$R^{*n}/s(R^{*n},q_{n})^r\defor A$
is clearly the $n$-fold coproduct of $R/s(R,q_{1})^r\defor A$ in the
category of those deformations $p:C\defor A$ which satisfy $s(C,\ker p)^r=0$.
From this latter fact, the definition of $h/B$,  the definition of 
the Grothendieck construction and the fact that the latter is the
homotopy colimit in CAT cf. [T1,1.2],
it follows immediately
that $h/B=\hoco_{\Delta^{op}}(n\mapsto S^{n+1})$. Here 
$S=\hom_{inf(C/A)}(R/s(R,q_0)^r,B)$ --which is a nonempty set because $R$ is cofibrant-- 
is thought of as a discrete category.
Hence taking nerves, $N(h/B)\approx\hoco_{\Delta^{op}}(n\mapsto NS^{n+1})$
(by [T1, 1.2]). But as $S$ is discrete this last $\hoco$ is nothing
but the simplicial set $n\mapsto S^{n+1}$ which is null homotopic.\qed
\enddemo
\bigskip
\proclaim{Corollary 5.2}Let $C$ be a category with deformations having
sufficient cofibrant objects, and $X:C@>>>\ssets$ a functor. Assume
$X(A)$ is fibrant for all $A\in C$. Given $A\in C$ choose a fibration
$f:RA\fib A$ with $RA$ cofibrant, and consider
$H_{inf}(-,X):A\mapsto \holi_{\Delta\times\Bbb N^{op}}X\check{f}$. Then
$H_{inf}(-,X)$ is a functor $C@>>>\ho\ssets$, is independent of the choices
made in its definition, maps deformations to isomorphisms and is equipped
with a natural map  $H_{inf}(-,X)@>>>X$. Further, if $A\in C'\subset C$ is 
a small subcategory of
interest containing a fibration $R\fib A$ with $R$ cofibrant, then
$H_{inf}(A,X)$ is the sheaf cohomology of $X$ on the infinitesimal
site $inf(C'/A)$ as defined in {\rm 5.0} above.
\endproclaim
\demo{Proof}Straightforward from the cofinality theorem for $\holi$
([BK,XI 9.2]) and the discussion 5.0.\qed
\enddemo
\bigskip
\subsubhead{5.3 Spectral Sequence }\endsubsubhead
Let $C$ be as in 5.2 above, and let $X$ be a sheaf of fibrant ssets.
We assume $\pi_0X,\pi_1X$ are abelian groups.
Fix a fibration $f:R\fib A$ as
in 5.2 and form the homotopy limit
$\holi_{\Delta\times\Bbb N^{op}}X\check{f}$. By [BK XI 4.3], the latter is
isomorphic to 
$\holi_{\Delta}(n\mapsto\holi_{\Bbb N^{op}}(m\mapsto X(R^{*n+1}/s(R^{*n+1},q_{n+1})^m)$. 
Thus by [BK, X 7.2] we have a  spectral sequence
$$
E_2^{r,s}=\pi^{r}\pi_{-s}\holi_{\Bbb N^{op}}(m\mapsto
X\check{f}(-,m))\tag{18}
$$
for $0\le r\le -s$, which, if $E_2^{r,-r}=0$, converges 
conditionally to
$H^{r-s}_{inf}(A,X)$. If $X$ is a sheaf of fibrant spectra then the
fringe constraints dissappear, and the sequence is always conditionally
convergent. A closely related spectral sequence is obtained as follows.
Write $H_{inf}(A,X)=\holi_{inf(C'/A)}X$ where $C'\subset C$ is any small
subcategory as in 5.2 above. Then by [BK, XI 7.1], we also have
a spectral sequence:
$$
{E'}^{rs}_2=H^r_{inf}(A,\pi_{-s}X)\tag{18'}
$$
with properties similar to those of \thetag{18}.
This may be regarded as a trivial kind of
Atiyah-Hirzebruch-Brown-Gersten ([BO]) spectral sequence. The sequences
\thetag{18} and \thetag{18'} agree in some cases; this is discussed
in 5.4.1 below. 
\bigskip
\proclaim{Lemma 5.4} The groups $E^{rs}_2$ are independent of the choice
of $f$.
\endproclaim
\demo{Proof} Let $f:R\fib A$ be a fibration with $R$ cofibrant. Let
$f\in C'\subset C$ be a small subcategory with deformations. Consider
the category $D(A)$ having as objects the fibrations $\fib A\in C'$
and as maps $(\alpha_0:B_0\fib A)@>>> (\alpha_1:B_1\fib A)$ the
pro-maps $B_0/s(B_0,\ker\alpha_0)^\infty@>>>B_1/s(B_1,\ker\alpha_1)^\infty$.
Give $D(A)^{op}$ the indiscrete topology, and consider 
$\alpha\mapsto\pi_s\holi_{\Bbb N^{op}}X(B/s(B,\ker\alpha)^\infty)$ as
a sheaf $\pi_s$ of abelian groups on $D(A)^{op}$. Then, by 2.1, we have 
$E_2^{0,s}=H^0(D(A),\pi_{-s})$. On the other hand the proof of 
[A, 3.1] shows that for $r\le 1$ the groups $E_2^{r*}$ vanish on injectives.
Summing up, $E_2^{r,s}=H^{-r}(D(A), \pi_{-s})$, and the lemma follows.\qed
\enddemo
\bigskip
\remark{Remark 5.4.1} If the shrinking functor is idempotent, then
the argument of the proof of the lemma above shows that 
$E_2^{r,s}={E'}_2^{r,s}$. In general if we have a
deformation
category structure in pro-$C$ extending that of $C$, (as in 1.7),
then, again by the proof of the lemma, 
$E_2^{r,s}=H^{r}_{pro-inf}(A,\pi_{-s}X)$, 
the cohomology of the sheaf $B@>>>\pi_s\holi_I XB_i$ on 
$inf(\PC/A)$. If in addition
the functor $C@>>>$Abelian Groups, $A\mapsto \pi_{-s}X(A)$ maps
deformations
into surjections, then 
$H^*_{pro-inf}(A,\pi_{-s}X)\cong H^*_{inf}(A,\pi_{-s} X)$, whence
 $E^2_{*,s}={E'}^2_{*,s}$. In
fact in
general if $G:C@>>>Ab$ is any functor and $f:R\fib A$ is a fibration
with cofibrant source, then the complex associated with the \v Cech 
pro-covering
$\check{f}$:
$$
\lim G(R/s(R,q_1)^n)@>\partial>>\lim G(R^{*2}/s(R^{*2},q_2)^n)
@>\partial>>\lim G(R^*3/s(R^{*3},q_3)^n)@>\partial>>\dots\tag"$C(A,G):$"
$$
computes $H^*_{inf}(A,G)$. This can be seen e.g. by mimicking the proof
that the \v Cech complex computes presheaf cohomology [A, 3.1]. The
complex $C$ can also be
used to compute infinitesimal
hypercohomology of some chain complexes. Indeed if
 $G:C@>>>$ ((positively graded chain complexes of
abelian groups)) is a functor which maps deformations into surjections,
then we have a weak equivalence $M(H_{inf}(A,SG))@>\sim>>STot_*(C(A,G))$.
Here  $C(A,G)$ is regarded as a (second quadrant) double complex, $Tot$
is the total chain complex, 
$S:$Chain Complexes$@>>>$S. Abelian Groups is the Dold-Kan functor
considered in 4.2 above and $M$ is its left adjoint, the (Moore) normalized
chain complex.
This fact follows from a long chain 
of homotopy equivalences which essentially use that $S$, having
a left adjoint (namely $M$), preserves limits, that it maps surjections
to fibrations, that $MS=1$, and that $\holi=\lim$
for both towers of fibrations (by [BK, XI 4.1-v)]) and cosimplicial abelian 
groups (by combining [BK, XI 4.4, X 4.9, 4.3, and 5.2-ii)]).
\endremark
\bigskip
\proclaim{Proposition 5.5} Let $C$ be a category with deformations
with idempotent shrinking functor. Let $X:C@>>>${\rm Fibrant }$\ssets$
be a functor. Assume the $\pi_nX$ are abelian groups for all $n$. Write 
$\check{X}(A)=X\check{1_A}$ for the composite
of $X$ with the \v Cech pro-covering {\rm 5.1} induced by the identity map. 
Then 
$$
H_{inf}(A,X)=\tot\bold{\check{X}}
$$
as homotopy types. If in addition, $\tot\bold{X}$ exists and $\tot\bold{X}\gamma(R)=\bold{X}(R)$
for cofibrant objects $R$, then $H_{inf}(A,X)=\tot\bold{X}(A)$.
The latter occurs iff the natural map $H^0_{inf}(R,\pi_nX)@>\cong>>\pi_nXR$ 
is an isomorphism for cofibrant $R$ and $n\ge 0$. In such case the 
$E_2$-term of the spectral sequence \thetag{18} vanishes for $r\ne 0$.
\endproclaim
\demo{Proof} The first assertion follows from  5.2 and 3.1-iv). If $\tot\bold{X}$
exists and agrees with $\bold X$ on cofibrants, then 
the cosimplicial group $n\mapsto\pi_r(R^{*n+1}/s(R^{*n+1},q_{n+1}))$ is
constant for cofibrant $R$. Hence the spectral sequence vanishes outside the
$y$-axis. It follows that the natural map $H_{inf}(R,X)@>>>XR$ is a
weak equivalence, and $H_{inf}(A,X)\cong\tot\bold{X}(A)$. The remaining assertion
follows from the fact that, by 5.4.1, $E'_2=E_2$ in this case.\qed
\enddemo
\bigskip
\remark{Remark 5.6}(Stratifying site and homotopy) The proposition above suggests an interpretation
of $H_{inf}$ in terms of homotopy, as defined in section 1. Indeed
we may think of $\check{X}$ as a homotopization of $X$. The construction
$X\mapsto\check{X}$ transforms any functor $X$ into one which maps deformation
retractions into weak equivalences. Hence in making the construction $\check{X}$
we are forcing $X$ into a functor which satisfies the conditions of 3.1. 
Back to the sheaf theoretic approach we may interpret $\check{X}(A)$ as
the sheaf hypercohomology on the stratifying site $strat(X/A)$, i.e. the site
consisting of those deformations $\defor A$ which are split. Note that $A$ is
cofibrant iff $strat(C/A)=inf(C/A)$. In terms of the stratifying site,
the proposition above states that $H_{inf}(A,-)=H_{strat}(R/s(R,I)^\infty,-)$
for any given fibration $R\fib A$ with cofibrant source having kernel $I$.
\endremark
\bigskip
\example{Application 5.7} (Tautological Character) Let $C$ be a category
with deformations having sufficient cofibrants. Assume for simplicity that
there is a functorial choice of fibration $\pi^A:RA\fib A$ with $RA$ cofibrant,
so that, given a functor $X:C@>>>$Fibrant Spectra, we can regard 
$H_{inf}(A,X)=\holi_{\Delta\times\Bbb N^{op}}X\check{\pi}^A$
as a functor to fibrant spectra, rather than just homotopy types. Here by
a fibrant spectrum  we mean a sequence
$X=\{ X^n:n\ge 0\}$
of fibrant ssets and weak equivalences $X^n\tilde{@>>>}\Omega X^{n+1}$.  
Write $\tau X^n:=(\hofiber (H_{inf}(A,X^{n+1})@>>>X^{n+1}))$. Then
we have a map  $c^\tau:X@>>>\tau X$; we call this the tautological
character. By definition it induces a weak equivalence
$X(B\defor A):=\hofiber (XB@>>>XA)@>\approx>>\tau X(B\defor A)$
of the relative spaces of a deformation. Hence $c^\tau$ satisfies
a tautological Goodwillie theorem \thetag{2}. Of course the identity
$X=X$ has the same property; unlike the identity however, the map $c^\tau$
is universal in the following sense. If $c:X@>>>Y$ is another map (character)
for which the Goodwillie theorem holds, then from the fact that $\holi$
preserves fibration sequences, it follows that the induced map
$\tau X@>\approx >>\tau Y$ is a weak equivalence. Hence $c^\tau$ factors
through $c$. 
In other words the tautological character is a coarser invariant 
than any other character with a Goodwillie theorem. 
\endexample
\bigskip

\subhead{6. Infinitesimal $K$-theory }\endsubhead
\bigskip
The purpose of this section is to apply the infinitesimal hypercohomology
machine to rational $K$-theory. First of all, choose any connective
functorial spectrum $K(A)$ with  homotopy $\pi_nK(A)=K_n(A)$,
$n\ge 0$. Then set $K^\rat=\rat_\infty K$. Most of the spectral sequence
\thetag{18} can be computed immediately using 4.3 and 5.5; we have:
$$
E_2^{rs}=\left\{\aligned LK^\rat_{-s}(A) &\quad r=0,s\le -2\\
K_0^\rat(A) &\quad r=s=0\\
0&\quad r\ne 0, s\ne -1\endaligned\right.
\tag{19}
$$
It follows that
$$
H^{-n}_{inf}(A,K^\rat))=LK^\rat_n(A)\qquad n\ge 2
$$
Note also that, with the possible
exception of the map $d_2:K^\rat_0(A)@>>>E^{2,-1}_2$, all spectral
differentials
are zero. Hence, for a full computation it suffices to compute this map and
the terms $E_2^{r1}$ $r\ge 0$.
By 5.4.1, we have:
$$
E^{r,-1}_2=H^{r}_{inf}(A,K^\rat_1)\tag{20}
$$
Next we use Goodwillie's theorem \thetag{16} to relate the groups 
\thetag{20} with
the corresponding groups for negative cyclic homology. We have
a commutative diagram with exact rows:
$$
\CD
0\rightarrow D(A,K^\rat_1)@>>>C(A,K^\rat_1)@>>>K^\rat_1(UA)\rightarrow 0\\
@VVV @VVV @VVV\\
0\rightarrow
D(A,HN^\rat_1)@>>>C(A,HN^\rat_1)@>>>HN^\rat_1(UA)\rightarrow0\\
\endCD\tag{22}
$$
Here the terms on the right hand side are the cochain complexes associated
to the constant cosimplicial objects, and the $D$-terms are the kernels
of the projections. The vertical maps are induced by the character of
Jones-Goodwillie. 
From 5.4.1, 5.6, and the long cohomology sequence of (22), 
we get the isomorphism:
$$
H^{r}_{inf}(A,K^\rat_1)\cong H^{r}_{inf}(A,HN^\rat_1) \qquad r\ge 2\tag{23}
$$
\bigskip
{\bf Conjecture 
6.0:} $H^n_{inf}(A,HN^\rat_1)=0$ for $n\ge 1$.
\bigskip
The conjecture is not stated for $r=0$ because in this case we prove
it below, so it is not a conjecture. We shall justify the conjecture
by showing that several things one expects should happen depend on
its validity. For instance note 
 that the groups \thetag{23} determine
the negative homotopy
of both
$H_{inf}(A,K^\rat)$ and
$H_{inf}(A,N^\rat)$,
where the latter is the hypercohomology of the connective spectrum associated
by Dold-Kan to the rational negative cyclic complex truncated below zero.
Hence the conjecture implies that this negative homotopy groups are zero.
 Moreover we have a commutative diagram:
$$
\CD
H(A,K^\rat)@>>>K^\rat(A)@>c_\tau>>\tau K^\rat(A)\\
@VH(-,ch)VV    @VVchV               @VV\tau chV\\
H(A,N^\rat)@>>>N^\rat (A)@>>>\tau N^\rat (A)\\
\endCD
\tag{24}
$$

Here  $c_\tau$ is the tautological character, the columns are induced by the
character of  Jones-Goodwillie and the last
of these is a weak equivalence by \thetag{16} and 5.7. Hence for $ch$ to
be the
tautological character we need  $H_{inf}^{*n}(A,HN_1)=0$ for all $n$, i.e.
we need the conjecture to hold for $A$.  
We prove next (Lemma 6.1) that $H^0_{inf}(A,HN_1)=0$. In theorem 6.2 below
we
use this vanishing result to 
extend the long
exact sequence  of theorem 4.3 so as to include $K_1$.
\bigskip

\proclaim{Lemma 6.1} $H^0_{inf}(A,HN^\rat_1)=0$.
\endproclaim
\demo{Proof} It suffices to prove the lemma for $\rat$-algebras.
 Consider the functor $\Omega^1_\natural$,
$A\mapsto\Omega^1A/[A,\Omega^1A]$.
I claim it suffices to show that $H^0_{strat}(A,\Omega^1_\natural)=0$ for
all $\rat$-algebras $A$. To prove the claim, proceed as follows. 
Note that, for a quasi-free pro-algebra,  the pro-cyclic mixed complex is
homotopy equivalent to the pro-de Rham mixed complex 
$MX=(\Omega^0\oplus\Omega^1_\natural, b,d_\natural)$. Hence 
$H^0_{inf}(A,HN^\rat_1)\cong\pi_1H_{inf}(A,NX)$, the hypercohomology
of the negative cyclic complex of $MX$. As $NX_1=\Omega^1_\natural$
and is zero in higher degrees, it follows that
$\pi_1H_{inf}(A,NX)\subset H^0_{inf}(A,\Omega^1_\natural)$
Now by 3.7, 5.4.1, and 5.6, we have 
$H^0_{inf}(A,\Omega^1_\natural)\cong H^0_{strat}(TA/JA^\infty,\Omega^1_\natural)$.
The claim follows from the fact that $\lim$ preserves kernels.
To prove $H^0(A,\Omega^1_\natural)=0$, we must show that if $\omega\in\Omega^1A$
is a $1$-form such that 
 the class of $(\partial_0-\partial_1)\omega$ in $\Omega^1(QA/qA^n)_\natural$
is zero for all $n\ge 0$, then $\natural\omega=0\in\Omega^1A_\natural$. We shall
show something stronger, namely if $(\partial_0-\partial_1)\omega$ is
zero in $\Omega^1(QA/qA^2)_\natural$ then $w_\natural=0$. We have
$QA/qA^2=A\oplus\Omega^1A$ and 
$\Omega^1(A\oplus\Omega^1A)=\Omega^1A\oplus\tilde{A}\otimes_\rat\Omega^1A
\oplus\Omega^1\otimes_\rat A\oplus\Omega^1A\otimes_\rat\Omega^1A$. 
We get $\Omega^1(A\oplus\Omega^1A)_\natural\cong\Omega^1A_\natural\oplus
\frac{\tilde{A}\otimes_\rat\Omega^1A}{M}\oplus\Lambda^2_{\tilde{A}}\Omega^1A$.
Here $M$ is the subspace generated by the elements of the form
$1\otimes b\omega a-b\otimes\omega a-a\otimes b\omega+ ab\otimes\omega$
and of the form $1\otimes [b,\omega]$. The relation which permits
eliminating the factor $\Omega^1A\otimes A$ is
$$
\omega\otimes a=1\otimes aw-a\otimes\omega\tag{25}
$$
It follows that the left multiplication map 
$\tilde{A}\otimes_\rat\Omega^1A@>>>\Omega^1A$
induces a projection 
$\Omega^1(A\oplus\Omega^1A)_\natural\fib\Omega^1A_\natural$.
Let $\Omega^1A\ni\omega=\sum a_idb_i$  ($a_i\in\tilde{A}$, $b_i\in A$).
Then $p(\partial_0-\partial_1)\omega=
p(\sum da_i\otimes b_i+a_i\otimes db_i+ da_i\otimes db_i)=\omega_\natural$
by \thetag{25}. \qed
\enddemo
\bigskip
\proclaim{Theorem 6.2} Let $A$ be a ring. 
Then:
\smallskip
\item{i)} The tautological character
$c^\tau_n:K^\rat_n(A)@>>>\tau K^\rat_n(A)$ coincides with the Jones-Goodwillie 
character $ch_n:K_n^\rat(A)@>>>HN_n^\rat(A)$ for $n\ge 2$. For 
$n=1$, the map $c^\tau_1$ factors as $ch_1$ followed by an injection.
\smallskip
\item{ii)} There is a natural map of spectra 
\noindent $f:\hofiber(ch:K^\rat (A)@>>>N^\rat (A))@>>>H_{inf}(A,K^\rat)$
such that the induced map $\pi_n(f)$ is an
isomorphism for $n\ge 1$.
\smallskip
\item{iii)} For $n\ge 2$ the infinitesimal hypercohomology groups agree
with the derived functor groups of {\rm Theorem 4.3}; we have
$H_{inf}^{-n}(A,K^\rat)=LK^\rat_n(A)$.
The long exact sequence {\rm 4.3-i)} extends to the right as follows:
$$
K^\rat_2(A)@>>>HN_2(A)@>>>H_{inf}^0(A,K^\rat_1)@>>>K^\rat_1(A)@>>>HN_1(A)
@>>>H^1(A,K^\rat_1)
$$
If furthermore, conjecture {\rm 6.0} holds for $A$, then the last map
above is surjective, and $H_{inf}^n(A,K)=\tau_nK=0$ for $n<0$.
\endproclaim
\demo{Proof} It follows from the lemma above and the exact homotopy sequence of the
commutative diagram of fibrations \thetag{24}.\qed
\enddemo
\bigskip
\subhead{Infinitesimal v. De Rham cohomology 6.2}\endsubhead
A theorem of Grothendieck [Dix, Th. 4.1] establishes that for a 
smooth scheme of characteristic zero, its de Rham cohomology is the
same as the infinitesimal cohomology of the structure sheaf:
$$
H^*_{dR}(X)\cong H^*_{inf}(X,O_X)\tag{26}
$$
We note that, even for a smooth affine scheme $Spec A/Spec\rat$, the 
infinitesimal
site in the sense of Grothendieck is larger than the infinitesimal
site of $A$ as an object of the category with deformations of 
commutative $\rat$-algebras. However the resulting cohomologies
are the same, cf. [Dix, isomorphism (*) on page 338]. 
Next we investigate a non commutative analogue of Grothendieck's
theorem. The non-commutative analogue of de Rham cohomology 
we use is the cohomology of the complex:
$$
DR:\quad HH_0(A)@>Bi>>HH_1(A)@>Bi>>HH_2(A)@>Bi>>\dots\tag{27}
$$
\bigskip
\proclaim{Theorem 6.3} (Compare \thetag{26}) Let $A$ be a quasi-free 
$\rat$-algebra. Then there is a natural map 
$$
H^n_{DR}(A)@>>>H^n_{inf}(A,O/[O,O])
$$
(Here we write $O/[O,O]$ for $HH_0$ to emphazise the relation of this
stament with Grothendieck's theorem \thetag{26}.) This map is an isomorphism
for $n=0$. If conjecture {\rm 6.0} holds for $A$ then it is an isomorphism
for all $n$.
\endproclaim
\demo{Proof} Consider the hypercohomology of the complex \thetag{27}.
As $A$ is quasi-free, this can be calculated from the \v Cech
pro-covering associated with the identity $A=A$. Thus $H_{inf}(A,DR)$
is the cohomology of the double cochain complex with as $n-th$ column
the $DR$-complex of the quasi-free pro-algebra $P^n=A^{*n+1}/q_{n+1}^\infty$.
Here $HH_r(P^n)=H_r(\lim_m C(P^n_m))$ is the homology of the limit
of the Hochschild complexes. As $P^n$ is quasi-free, we have $HH_r(P^n)=0$
for $r\ge 2$. On the other hand, by the tubular neighborhood theorem 
([CQ1, Th2])  $P^n\cong T_{\tilde{A}}(M_n)/<M_n>^\infty$, where 
$M_n=\Omega^1A\oplus\dots\oplus\Omega^1A$
 ($n$ factors). It is not hard to see that in general, if $B=B_0\oplus B_1\oplus\dots$
is a graded algebra, then for each $r$ the pro-vectorspace
$\{HH_r(B/B_+^m):m\in\Bbb N\}$ is isomorphic to the pro-vectorspace
$\{\oplus_{i=0}^m HH^{deg=i}_r(B):m\in \Bbb N\}$ which is clearly
Mittag-Leffler. In our case, this implies that 
$HH_*(P^n)=\lim_m HH_*(P^n_m)$. Hence $H_{inf}^*(A,DR)$ is the 
cohomology of the double cochain complex:
$$
\CD
HH_1(A)@>>>\lim_m HH_1(P^1_m)@>>>\lim_m HH_1(P^2_m)@>>>\\
@AdAA              @AdAA              @AdAA \\
HH_0(A)@>>>\lim_m HH_0(P^1_m)@>>>\lim_m HH_0(P^2_m)@>>>\\
\endCD\tag{28}
$$
By 5.4.1, the rows of this complex compute $H_{inf}^*(A,HH_i)$.
Hence we have a spectral sequence:
$$
{E"}_1^{r,s}=H_{inf}^r(A,HH_s)\implies H^{r+s}_{inf}(A,DR)\tag{29}
$$
On the other hand if we first take cohomology of the columns,
and then of the rows, we get $HP_i(A)$ in the $(0,i)$ entry
and zero elsewhere. Hence $H^n_{inf}(A,DR)=HP_n(A)$ for $n=0,1$
and is zero if $n\ge 2$. The theorem now
follows from the convergence of the spectral sequence \thetag{29},
and the fact that
$\operatornamewithlimits{lim}_m HH_1(P^n_m)=HH_1(P^n)=HN_1(P^n)=
\operatornamewithlimits{lim}_mHN_1(P^n_m)$.\qed
\enddemo
\bigskip

\define\CAT{\operatorname{CAT}}
\define\heq{\overset\sim\to\rightarrow}
\define\HOM{\operatorname{HOM}}

\subhead{7. Categorical character}\endsubhead

Throughout this section
$C$ is a fixed small category with deformations having sufficient
cofibrant objects. We consider functors from $C$ to 
$\CAT$, the large category of all small categories.
\bigskip
\subhead{Homotopy limits and colimits in CAT 7.0}\endsubhead
By the {\it homotopy colimit} of a functor $F:I@>>>\CAT$ ($I\in\CAT$)
we mean the Grothendieck construction:
$$
\hoco_I F:=\int_IF
$$
and by the {\it homotopy limit}
we mean the pullback:
$$
\CD 
\holi_IF@>>>0\\
@VVV @VVidV\\
\HOM(I,\hoco_IF)@>>\pi_*>\HOM (I,I)\\
\endCD
$$ 
Here $\HOM$ is the
functor category, $\pi_*$ is induced by the natural projection
$\pi:\hoco_IF@>>>I$, $0$ is the category with only 
one map and $id$
maps the only object of $0$ to the identity functor.
In other words,
$\holi_IC$ is the category whose objects are the functors
$s:I@>>>\hoco_IC$ such that $\pi s=1$ and whose maps are the natural
transformations which project to identity maps through $\pi$.
Thomason ([T1]) showed that, upon taking nerves, the $\hoco$ defined above
has the same homotopy type as its simplicial counterpart,
i.e. $N\hoco_IF\approx\hoco_INF$. On the other hand it was observed
in [L, p74] that there is an isomorphism of simplicial sets
 $N\holi_IF\cong\holi_INF$. By definition a map of categories
is a weak equivalence if its nerve is a weak equivalence of SSets.
Hence by [BK,XII 4.2], $\hoco$ preserves weak equivalences. Similarly
$\holi$ preserves weak equivalences between categories having fibrant
nerve, by [BK, XI 5.6 ]. Moreover it is proven in [L, Th. 1] that $\holi$
maps adjoint functors to adjoint functors, and thus to weak equivalences.
Here is a description of both $\hoco_IF$ and $\holi_IF$ in terms of
objects
and arrows. An object of $\hoco_IF$ is a pair $x_i:=(i,x)$ where $i\in I$ and
$x\in F_i=:F(i)$. A map $x_i@>>>y_j$ is a pair 
$(\alpha,\rho)$
with $\alpha:i@>>>j\in I$ and  $\rho:\alpha x@>>>y\in F_j$. Here
we abbreviate $F(\alpha)(\rho)(\mu)$ as $\alpha(\rho)\mu$. 
Hereafter we shall omit $F$'s and $(,)$'s whenever no 
confusion is possible. Composition is defined as in a semidirect 
product: $(\alpha,\rho)(\beta,\mu)=(\alpha\beta, \rho\alpha(\mu))$.
On the other hand, $\holi_IC$ is the category of
all pairs of families 
$(x,\rho):=(\{ x_i\}_{i\in I}, \{\rho_\alpha\}_{\alpha\in I})$, indexed 
respectively
 by the objects
and the maps of $I$, where $x_i\in F_i$, and if
$\alpha:i@>>>j$ is a
map in $I$, then $\rho_\alpha:\alpha x_i@>>>x_j$
is a map in
$F_j$. The family of $\rho$'s is subject to the conditions: 
$$
\rho_1=1
\text{\qquad\qquad}\rho_{\alpha\beta}=\rho_\alpha\alpha(\rho_\beta)
 $$
In the first
identity, the $1$ on the left is an identity map $1:i@>>>i$
 and the $1$ on the right is the identity of $x_i$ in $F_i$; in the second
identity, $i_0@<\alpha<<i_1@<\beta<<i_2$ are composable maps in $I$
and the identity is of maps $\alpha\beta x_{i_2}@>>>x_{i_0}$ in
$F_{i_0}$. A map $f:(x,\rho^x)@>>>(y,\rho^y)$ in $\holi_IF$ is a
family of maps $f_i:x_i@>>>y_i\in F_i$ indexed by the objects of
$I$ such that the following diagram commutes for every map
$\alpha:i@>>>j\in I$:
$$
\CD
\alpha x_i@>f_i>>\alpha y_i\\
@V\rho^xVV @VV\rho^yV\\
x_j@>>f_j>y_j\\
\endCD
$$
\definition{ Sheaves of objects of a functor $F:C@>>>\CAT$ 7.1} Let $C$ be a 
small category
with deformations and let $F:C@>>>\CAT$ be a functor. Write:
$$
F_{inf}(A):=\holi_{inf(C/A)}F
$$
We call $F_{inf}(A)$ the category of {\it sheaves of objects} of $F$ on
$inf(C/A)$. If $F$ is constant, we recover the usual notion of sheaves
as functors $inf(C/A)@>>>F(A)$. Because $\holi$ commutes with nerves, 
we have a map:
$$
NF_{inf}(A)@>>>H_{inf}(A,NF)
$$
which is natural up to homotopy. This map is a weak equivalence if
$NF$ is fibrant, i.e. if $F$ takes values in the category of small
groupoids=categories where every map is an isomorphism, but not in
general.
We note that as $\holi$ preserves adjointness,
the functor $F_{inf}$ maps deformations into weak equivalences. Moreover
we have a natural evaluation map $\epsilon:F_{inf}(A)@>>>F(A)$, 
$(x,\rho^x)\mapsto x_A$, $\tau\mapsto\tau_A$. It is straighfoward to
check that ${F_{inf}}_{inf}\cong F_{inf}$, whence $\epsilon$ is (uni)versal
among natural maps $G@>>>F$ where $G_{inf}\cong G$. In other words
the construction $F\mapsto F_{inf}$ has properties similar to those
of $H_{inf}(-,-)$.
\enddefinition
\bigskip
\example{Example 7.2} (Stratified objects) Next
we give an explicit description of $F_{inf}(A)$ for a groupoid functor
$F$ and cofibrant $A$, up to
category equivalence. In the case when $C$ is the category of commutative
algebras of finite type over a field, and $F(A)$ is the isomorphism
category of finitely generated modules we recover [BO, Prop. 2.11] 
(see also [Dix, 4.2]). 
Write
$\partial^n_i:A@>>>P^1_n=A*A/s(A,qA)^n$ for the $i$-th coface map
($i=0,1$). Given an object $(x,\rho^x)\in F_{inf}$, consider the
family of isomorphisms: 
$$
\epsilon_n:\partial_0^nx\heq\partial^n_1x\in F(P^1_n)
$$
defined by $\epsilon_n=\rho^{-1}_{\partial_1}\rho_{\partial_0}$.
Then:
\item{i)} $\epsilon_0=1$ 
\smallskip
\item{ii)} The family $\{\epsilon_n:n\ge 0\}$ is 
compatible with the maps $P^1_{n+1}\defor P^1_n$.
\smallskip
\item{iii)}(Cocycle Condition) The following identity holds:
$$
\partial_2(\epsilon_n)\partial_0(\epsilon_n)=\partial_1(\epsilon_n)
$$
We call the family of $\epsilon$ a {\it stratification} on $x$.
Conversely, an object $x\in FA$ with a stratification $\epsilon_*$,
yields an object of $F_{inf}(A)$ as follows. Given $B\defor A\in 
inf(C/A)=strat(C/A)$
choose a section $s_{B}$ and set $x_{B}=s_{B}x$. If $\alpha:B@>>>C\in inf(C/A)$
is a map, let $h:P^1_n@>>>C$ be a homotopy $\alpha s_{B}\equiv s_{C}$,
and set $\rho_{\alpha}=h(\epsilon_n):\alpha x_{B}@>>>x_{C}$. It is straightfoward
to check that these assignments define mutually inverse equivalences 
between $F_{inf}A$ and the category of objects of $FA$ equipped with 
a stratification. The same argument shows that, even if $A$ is not
cofibrant --but $F$ is still a groupoid-- then the category of stratified
objects is equivalent to $F_{strat}=\holi_{strat(C/A)}F$ (compare [Dix, 4.2]).
\endexample
\bigskip
\subsubhead{Categorical Character 7.3}\endsubsubhead
Let $M$ be a permutative monoidal category, and consider the simplicial
category $M^+:n\mapsto M^{n+2}$ of [T1, 4.3.1]. Form the
fibrant spectrum $n\mapsto Sp_nM$ of simplicial sets associated to the topological
spectrum constructed in [T1, 4.2.1]; the $0$-th space of this
fibrant spectrum has the weak homotopy type of the
nerve of $\hoco_{\Delta^{op}}M^+$, which is a categorical model
for the group completion of the realization $BM=|NM|$. Set:
$$
KM:=Sp_0M
$$
Now suppose
$f:M@>>>N$ is a functor of permutative monoidal categories, preserving
products in the strong sense, i.e.  up to natural isomorphism, and assume
$M$ is a groupoid. Then by 
[T2 5.2] there is another permutative category $P(f)$ and a functor
$N@>>>P(f)$ such that the sequence of base spaces
$Sp_0M@>>>Sp_0N@>>>Sp_0P(f)$ is a fibration up to homotopy. Use
this fibration sequence to define relative groups as follows. Let
$A\mapsto MA$ be a functor going from a category $C$ with deformations into
the category of (small) permutative monoidal categories which are groupoids,
with as maps the strong product preserving functors. 
Given $f:A@>>>A'\in C$, write:
$$
KM(f):=KP(f) \text{\quad and\quad} KM^{rel}A=KP(M_{inf}A@>>>MA)
$$
\bigskip
\proclaim{Tautology 7.4} Let $C$ be a category with deformations. Let 
$M:C@>>>\CAT$ be a functor, mapping objects to groupoids which are
permutative monoidal categories and maps to functors preserving products
in the strong sense, i.e. up to natural isomorphism. Then there is a 
natural character
$c:KM(A)@>>>KM^{rel}(A)$ which is induced by a map of categories,
has  $KM^{inf}A$ as fiber, and
is such that any deformation $f:B\defor A\in C$ induces an isomorphism 
of the relative groups $KM_n(f)\cong KM_n^{rel}(f)$ $n\ge 0$. This 
character fits
into a homotopy commutative diagram with homotopy fibration rows:
$$
\CD
KM_{inf}A@>>>KMA@>c>>K(M^{rel}A))\\
@VVV           @V1VV    @VVV\\
H_{inf}(A,KM)@>>>KMA@>c^\tau>>\tau KMA\\
\endCD
$$
Here the bottom row is the homotopy sequence of the tautological character
5.7.
\endproclaim
\bigskip
\remark{Remark 7.5} In the discussion above, we applied the $inf$ construction
{\it before} group completing; we considered 
$\hoco_{\Delta^{op}}M^+_{inf}$. The reader may wonder what happens if 
in place of the latter category we use
$(\hoco_{\Delta^{op}}M^+)_{inf}$. The fact is that the two categories
are homotopy equivalent. This is a problem of interchanging $\holi$ and
$\hoco$, which in general is not possible, but which in the particular
case where the index category of $\holi$ has a final object, yields
homotopy equivalent spaces (cf. [C2]). 
\endremark
\bigskip
\example{7.6 Case of rings: $K$-theory of connections} The tautology above 
applies notably
in the case when $C$ is a category of rings and $MA=\coprod Gl_nA$.
For $A$ cofibrant,  $(\coprod Gl_n)_{inf}A$ can
be described, upon the identification 7.2, as the category having
as objects the pairs $(m,n\mapsto\epsilon_n)$ where the first
coordinate is a non-negative integer and the second is a family
of matrices $\epsilon_n\in Gl_m(P^1_n)$ satisfying i)-iii) above.
There is no map $(m,\epsilon)@>>>(r,\theta)$ unless $r=n$; if $r=n$
a map is a matrix $\alpha\in Gl_m(A)$ with 
$\theta\partial_0(\alpha)=\partial_1(\alpha)\epsilon$.
In the case when $C$ is the category
of commutative algebras essentially of finite type over a field
of characteristic zero 
and $A\in C$ is smooth, a stratification on $A$ is the same
thing as a flat connection (cf. [BO]). Thus the category
$(\coprod Gl_n)_{inf}A$ is identified with the category having
as objects the pairs $(n,\nabla)$ where $\nabla$ is
a flat connection on $A^n$. There are no maps $(n,\nabla)@>>>(n',\nabla')$ 
if $n\ne n'$, and if $n=n'$, a map is a matrix
$\alpha\in Gl_n(A)$ such that
$\nabla'\alpha=\alpha\nabla$. We remark that isomorphism classes of 
connections in this sense
are the classical gauge equivalence classes [BL]. 
Next consider the case when $C$ is some small category of associative
but not necessarily commutative algebras over a characteristic zero
field, and let $A$ be quasi-free. As indicated above objects of
$(\coprod Gl_n)_{inf}A$ are stratified free modules $(m,\epsilon)$. 
Because $\epsilon_0$
is to be the identity, the map
$\epsilon_1:P^1_1=A\oplus\Omega^1A@>>>A\oplus\Omega$ is necessarily
of the form $\epsilon_1=1+\nabla$ where $\nabla$ is a right connection
on $A^m$. I tried to prove that the cocycle condition is equivalent
to this right connection being flat, but quit overwhelmed by
the horrendous calculations. I got as far as proving
that flatness is equivalent to extending $\nabla$ to a cocycle $\epsilon_2$.
The cocycle identity for $\epsilon_3$ already takes several pages to write
down. This seems to indicate one should use a deeper argument
than just brute force. In the commutative case such deeper 
argument comes from 
the interpretation
of stratified modules as $D$-modules ([BO 2.11-3)]). I have not been able to find
a good analogy of this interpretation in the non-commutative case. 
\endexample

\enddocument